\documentclass[reqno]{amsart}
\usepackage{comment}
\usepackage{mathrsfs}
\usepackage[dvipsnames]{xcolor}
\usepackage{mathtools}
\usepackage{centernot}
\usepackage[colorlinks,linkcolor=blue,citecolor=blue]{hyperref}
\usepackage{cleveref}
\allowdisplaybreaks
\usepackage{graphics}
\usepackage{tikz}
\usepackage{nicematrix}
\usepackage{bbm}
\usepackage{amsmath, amssymb}

\numberwithin{equation}{section}

\theoremstyle{plain}
\newtheorem{thm}{Theorem}[section] 
\newtheorem{lem}[thm]{Lemma}
\newtheorem{prop}[thm]{Proposition}

\theoremstyle{definition}
\newtheorem{Def}[thm]{Definition}

\theoremstyle{remark}

\newcommand{\ep}{\epsilon}
\newcommand{\gsl}{\fsl{g}}

\graphicspath{ {Plots/} }

\input Tdef2.def

\begin{document}

\title[Global existence for semi-linear hyperbolic equations.]{Global existence for semi-linear hyperbolic equations  in a neighbourhood of future null infinity}

\author[J.A. Olvera-Santamaria]{J. Arturo Olvera-Santamaria}
\address{Department of mathematics and statistics\\
University of Calgary, \\
2500 University Dr Nw, Alberta, Canada}
\email{arturo.olvera@ucalgary.ca }

\begin{abstract}
In this paper,  we establish the global existence of a semi-linear class of  hyperbolic equations in  $3+1$ dimensions,  that satisfy the \emph{ bounded weak null condition}.  We propose a conformal compactification of the future directed null-cone in  Minkowski spacetime,  enabling us to establish  the solution to the wave equation in a neighbourhood of  future null infinity.  Using this framework, we formulate  a conformal  symmetric hyperbolic Fuchsian  system of equations.   The existence of solutions to this Fuchsian system  follows from an application of  the existence theory developed in \cite{Oliynyk2021107},  and \cite{BOOS:2020}.  
\end{abstract}
\maketitle
\section{\bf Introduction}
Significant advancements have been made in the study of  the existence of solutions to hyperbolic equations,  from the 1960's to the present.   A notable breakthrough was achieved by Christodoulou \cite{Christodoulou:1986} and Klainerman \cite{Klainerman1986}, who proved the global existence of non-linear hyperbolic equations that satisfy the so-called \emph{null condition}.   Equations satisfying this condition exhibit global  solutions that decay like solutions of linear wave equations.   Although the results in \cite{Christodoulou:1986} and  \cite{Klainerman1986} were groundbreaking,   Choquet-Bruhat \cite{choquet2000null} proved that these results do not apply to the Einstein's equations and that there is no natural generalization of the null condition for them.   

Later,  Lindblad showed \cite{lindblad1992global} that there exist quasilinear  equations that do not satisfy the null condition but still admit global solutions that decay slower than solutions of linear wave equations.  In a subsequent paper,  Rodniansky and Lindblad \cite{lindblad2003weak} designed a more general condition, which  they called the \emph{weak null condition}, and demonstrated that is satisfied by the Einstein's equations.  Then they  used the weak null condition  to prove a global existence result for Einstein's equations in wave coordinates \cite{LindbladRonianski:2005}.   

The weak null condition is based on the idea that a certain class of  non-linear hyperbolic equations is asymptotically  governed  by an ODE.  Therefore,  if the solutions of the asymptotic ODE exist,  have initial data decaying sufficiently fast,  and grow at most exponentially,  then the original system also admits a global solution.  The null condition can be viewed as a specific case of the weak null condition.  It remains an open problem to determine  if all non-linear hyperbolic equations satisfying the weak null condition have global solutions.  

The weak null condition appears to be very general,  leading  authors  to focus on specific cases,  see for example \cite{lindblad2018global,lindblad2020local,lindblad2022weak,Alinhac:2003,Alinhac:2006,Bingbing_et_al:2015,DengPusateri:2018,HidanoYokoyama:2018,Katayama_et_al:2015,Keir:2019}.  In   \cite{Keir:2018},  Keir proved the global existence for solutions to quasilinear wave equations satisfying the weak null condition along with a  hierarchical structure in the semi-linear terms.  In  \cite{Oliynyk2021107},  we established the  global existence of semi-linear wave equations that satisfy a restricted version of the weak null condition,  which we call the \emph{bounded weak null condition}.  This version  includes Keir's hierarchical condition.  More importantly,  in \cite{Oliynyk2021107},  our initial data does not require to be compactly supported,  unlike in  \cite{Keir:2018}.  In this paper,  we complement the results from \cite{Oliynyk2021107} by proving the global existence  of semi-linear wave equations of the form
\begin{equation}  
\gb^{\alpha\beta}\nablab_\alpha \nablab_\beta \ub^K = \ab^{K\alpha\beta}_{IJ}\nablab_\alpha \ub^I \nablab_\beta \ub^J, \label{Mbwave} 
\end{equation}
on a neighbourhood of future null infinity,  where the  $\ub^I$,  are the components of the unknown $\ub$.  The  region of Minkowski  spacetime that we are interested in,  is  the future directed null-cone in $\mathbb{R}^4$ with origin at $\xb^i=0$,  that is
\begin{equation}\label{M+region}
\Mb= \{(\xb^i) \in \mathbb{R}^4 \ | \ \xb^0>0, \ \gb_{ij}\xb^i\xb^j<0 \}.  
  \end{equation}  
There exist a $N$ rank vector bundle $V$ such that the unknown $\ub$ is a  section of $V$,  and  $1\leq I,J,K \leq N$  \footnote{See Appendix \ref{indexing} for our indexing conventions.}. The  $\ab_{IJ}^{K}=\ab_{IJ}^{K\alpha\beta}\delb{\alpha}\otimes\delb{\beta}$,  are prescribed smooth (2,0)-tensors fields on $\Rbb^4$,  and $\nablab$ is the Levi-Civita connection of the Minkowski metric.   We use the notation $(\hat{x}^{\mu})$ to denote Cartesian coordinates,  and $(\xb^{\mu})$ to denote spherical coordinates
\begin{equation*}
(\xb^\mu) = (\xb^0,\xb^1,\xb^2,\xb^3)=(\tb,\rb,\thetab,\phib),
\end{equation*}
which we use to write the  Minkowski metric
\begin{equation} \label{gbdef}
\gb =  -d\tb\otimes d\tb + d\rb \otimes d\rb + \rb^2\gsl,
\end{equation}
where
\begin{equation}\label{gsldef}
\gsl = d\thetab \otimes d\thetab + \sin^2(\thetab) d\phib \otimes d\phib,
\end{equation}
is the canonical metric on the $2$-sphere $\mathbb{S}^2$. 
For simplicity\footnote{The results of this article can be generalized to allow non-covariantly constant  tensors $\ab_{IJ}^{K}$ provided that they satisfy suitable asymptotics.},  we assume that the tensor fields $\ab^K_{IJ}$ are covariantly constant,  i.e.  $\nablab \ab_{IJ}^{K}=0$,  which is equivalent to the condition that the components of $\ab_{IJ}^{K}$ in a Cartesian coordinate
system $(\xh^\mu)$ are constants, that is, $\ab_{IJ}^{K} = \ah_{IJ}^{K\alpha\beta}\hat{\partial}_{\alpha}\otimes\hat{\partial}_{\beta}$
for some set of constant coefficients $ \ah_{IJ}^{K\alpha\beta}$.  Moreover,  equation \eqref{Mbwave} satisfies the \emph{bounded weak null condition},  which means that  the solutions to the the asymptotic equation \eqref{asympprop1.1} defined below, exist and are bounded.     

We follow the techniques and structure outlined in  \cite{Oliynyk2021107},  to establish the global existence of solutions to \eqref{Mbwave} in a subset of \eqref{M+region}.  Most of this paper is devoted to write the equations \eqref{Mbwave} into a Fuchsian system of symmetric hyperbolic equations that satisfies suitable conditions to apply the theory developed in  \cite{Oliynyk2021107} and \cite{BOOS:2020}.  One of the key ingredients in this paper is the introduction of a compactification of the future directed null cone,  such that the theory developed in \cite{Oliynyk2021107} is applicable to a system on this space.   In Subsection \ref{FMet}, we briefly describe the Fuchsian method  and in Subsection \ref{BNWC}, we define the version of the bounded weak null condition used in this work.  In Section  \ref{CMMS},  we introduce a  conformal map that compactifies the outgoing null-rays, enabling us  to approach future null infinity in a neighbourhood of $\rb=0$.  Using  the conformal transformation \eqref{map},  we map  the region \eqref{M+region} on to $(0,\infty)\times (0,1)\times \mathbb{S}^2$ and push-forward  the wave equation \eqref{Mbwave} onto the  new conformal manifold,  which leads to the conformal wave equation \eqref{weq1},  and \eqref{eweq}  defined on the manifold  \eqref{M0}, whose closure is compact.  

Then,  we  define the variables $U^K=(U^K_\Ic)$  by equations  \eqref{chv1}, which transform the system \eqref{eweq} into a first order system.  We then   propose  the change of variables \eqref{chvm},  to ensure that the resulting first order system is symmetric hyperbolic.  This change of variables condenses the most singular terms into a  multiple of the semi-linear term $V^I_0V^J_0$.  Additionally,  it reveals  the form of the asymptotic equation \eqref{aeqi3},  which involves the most singular term.  However, this singular term  does not pose a significant challenge  since it is  removed later using the  flow of the asymptotic equation \eqref{aeqi2} to redefine $V_0$ to a new variable $Y$ determined implicitly  by the flow equations \eqref{asympIVP.1},  \eqref{asympIVP.2}, \eqref{Ydef}.  The removal of the most singular term leads to the evolution equation \eqref{MwaveM},  which becomes part of the complete Fuchsian system \eqref{MwaveO}-\eqref{Jcdef}.   

The extended system is defined by the equations  \eqref{MwaveE},  \eqref{FcKdef},  \eqref{QKdef},  \eqref{Gcdef},   \eqref{GcK0def},  \eqref{GcK1def},  \eqref{GcK2def}\eqref{PbbGc2=0},   on the closed manifold $(0, t_0) \times \mathcal{S} $ (see  section \ref{extendst}  eq. \eqref{extst}),  which is a key requirement  for the Fuchsian method \cite{BOOS:2020}.   In Subsection \ref{differsys},  we  differentiate the extended system \eqref{MwaveE} to derive the system \eqref{difs}, \eqref{Wvar},  \eqref{Qcdef},  \eqref{HcKdef}  which is taken as an evolution equation for the variables $W^K_{j}=t^{\kappa}(\Dc_j V^K)$.  These equations are also part of the complete Fuchsian system \eqref{MwaveO}-\eqref{Jcdef}.  Subsequently,   we  apply the projection operator $\Pbb$ to the extended system \eqref{MwaveE},  yielding an equation for the variable  $X^K=t^{-\nu}\Pbb V^K$.   Then,  we combine the three systems 
\eqref{difs},   \eqref{MwaveN}  \eqref{MwaveM} involving the variables $W^K_j, X^K,Y^K$ into the single Fuchsian system \eqref{MwaveO} to obtain an evolution system for $Z=(W^K_j, X^K,Y^K)$.   Finally,  in Section \ref{complete}, we show that under the flow assumptions from section \ref{Asymptoticflowassump},  the Fuchsian system \eqref{MwaveO} satisfies all the necessary conditions to apply the  Global Initial Value Problem (GIVP) existence theory from \cite{BOOS:2020}.  Applying Theorem 4.1 from \cite{BOOS:2020},  we establish the  GIVP result for the Fuchsian system  \eqref{MwaveJ.1}, \eqref{MwaveJ.2},  which  by construction,  implies a global existence result for the original system of wave equations \eqref{Mbwave}, for sufficiently small initial data.

\subsection{The Fuchsian Method} \label{FMet}
The results presented in this paper are part of a broader  research program that  employs the  Fuchsian method as a tool to prove the global existence of solutions to non-linear hyperbolic equations in various settings.  The essential idea of the method  is to transform a non-linear system of hyperbolic equations into  a Global Initial Value Problem (GIVP) for a  first order  Fuchsian system of  symmetric hyperbolic equations.  This is achieved  by  applying a suitable  conformal transformation to the original system of equations.  Then,  using  energy estimates,  we prove the global existence of solutions to the conformal  equations.  By construction,  these solutions  yield  the  global existence of solutions to the original set of equations.   Examples of  GIVP applications  can be  found in \cite{Oliynyk2021107,BOOS:2020, LeFlochWei:2020,LiuOliynyk:2018b,LiuOliynyk:2018a,LiuWei:2019,Oliynyk:2020,Wei:2018,FOW:2020,LiuWei:2019prp}.  This method is notable for it is \emph{simplicity} compared to other techniques,  and its capacity to handle singular systems of hyperbolic equations.   The GIVP offers significant advantages over the Singular Initial Value Problem (SIVP),  which requires to prescribe asymptotic data at the singular time.  In contrast,  in a GIVP  we prescribe initial data at a finite time $t=t_0$ with  the challenge being   to prove that  solutions to the system of wave equations exist up to the singular time.   This  makes it a promising method  to study  singular systems of  equations where initial data near the singularity is unknown.  Readers interested in the SIVP may consult for example  \cite{andersson2001,choquet2006half,choquet2005topologically,damour2002,heinzle2012,IsenbergKichenassamy:1999,isenberg2002,kichenassamy1998,OliynykKunzle:2002a,OliynykKunzle:2002b,Rendall:2004}.  

The first step in the Fuchsian method,  consists in transforming a system of non-linear hyperbolic equations,  for example eq.  \eqref{Mbwave} into a first order  symmetric hyperbolic system of the form 
\begin{equation}
B^0(t,u)\del{t}u + B^i(t,u)\nabla_{i} u  = \frac{1}{t} \Bc(t,u)\Pbb u + F(t,u), \hspace{1cm}  \text{in} (t_1,t_0]\times \Sigma,\label{symivp}
\end{equation}
where,  the unknown $u$ is a time-dependent section of  a N rank vector bundle $V$.  The matrices  $B^0,  B^i$ are symmetric operators on $V$.   $\Bc$ is a linear operator on $V$, and $\Pbb$ is a time-independent,  covariantly constant, symmetric projection operator.   Then the system  \eqref{symivp} is viewed as a global initial value problem  with suitable initial data specified at $ u_0 \in  {t_0}\times \Sigma$,  and the objective is to establish the existence of solutions to \eqref{symivp} in an interval that reaches the singular time at $t=0$,  that is $t\in (0,t_0]$.   It is important to note that the system \eqref{symivp} is defined in the non-phyiscal spacetime,  in other words, it is defined in the conformal version of the initial spacetime.   We  apply the existence  theory  for Fuchsian systems from \cite{Oliynyk2021107},  to the system \eqref{symivp}  provided that it  satisfies the structural conditions given in \cite{Oliynyk2021107},  and  \cite{BOOS:2020}.  Although the  transformation process is  particular to each system,   we can highlight  4 main steps required to  transform a wave equation into a Fuchsian system:
\begin{enumerate}
\item Transforming the physical manifold where the  original system is defined,  into a bounded N-dimensional non-physical manifold whose boundary represents infinity of the physical manifold.  In \cite{Oliynyk2021107},  we carried out this step  by applying   Friedrich's  \emph{cylinder at infinity}   conformal transformation \cite{Friedrich:JGP_1998}.  
\item Transforming the second order conformal wave equation into a first order symmetric  hyperbolic equation. 
\item A rescaling on time might be required in order to meet the coefficient assumptions from \cite{BOOS:2020}. 
\item  Verification of the structural conditions for a Fuchsian system in order to apply the Theorem 4. 1 from \cite{Oliynyk2021107}.
\end{enumerate}

\subsection{The bounded weak null condition}\label{BNWC}
To continue we introduce  the out-going  null one form $\Lb= -d\tb +d\rb$, and we define  the functions
\begin{equation}\label{nullbk}
\bb^K_{IJ}=\ab^{K}_{IJ}\Lb_{\mu}\Lb_{\nu}=\ab^{K00}_{IJ}-\ab^{K01}_{IJ}-\ab^{K10}_{IJ}+\ab^{K11}_{IJ}, 
 \end{equation} 
the importance of the term $\bb^K_{IJ}$ is that  the \emph{null condition} is satisfied when $\bb^K_{IJ}=0$.  With the help of  the functions \eqref{nullbk} we define the asymptotic equation associated to \eqref{Mbwave} by
\begin{equation}\label{aseQ}
\del{t}\xi=\frac{1}{t}Q(\xi)
\end{equation}
where 
\begin{equation}\label{aeqi}
Q(\xi)=(Q^K(\xi))\coloneqq  -\frac{\chi(\rho)\rho^3m}{2(\rho^{2m}-(1-\rho^m)^2)t}\bb^K_{IJ}\xi^I \xi^J. 
\end{equation}
The asymptotic equation \eqref{aseQ} is defined in terms of the coordinates $(t,\rho,\theta,\phi)$,  which arise from the compactification \eqref{map} of  a neighbourhood of future null infinity and the rescaling of $r$ \eqref{rtr}.  Our time coordinate $t$ is such that at  $t_0$ (see equation \eqref{time0}),   we set suitable initial data,  and the time  $t=0$ corresponds to the evolution of the initial data towards future null  infinity.  We say that the equation \eqref{Mbwave} satisfies  \emph{the bounded weak null condition}  if the solutions to the asymptotic equation \eqref{aeqi} exist and are bounded.
\begin{Def} \label{bwnc}
\noindent The asymptotic equation  satisfies the \textit{bounded weak null condition} if there exist constants $\Rc_0>0$ and $C>0$ such that solutions of the asymptotic initial value problem (IVP)
\begin{align}
\del{t}\xi &= \frac{1}{t}Q(\xi), \label{asympprop1.1} \\
\xi|_{t=t_0} &= \mathring{\xi}, \label{asympprop1.2}
\end{align}
exist for $t\in \left(0,t_0\right]$ and are bounded by $\dsp \sup_{0<t\leq t_0}|\xi(t)| \leq C$
for all initial data $\mathring{\xi}$ satisfying $|\mathring{\xi}| < \Rc_0$. 
\end{Def}

\section{\bf Conformal mapping of Minkowski spacetime  near future null infinity}\label{CMMS}
In  \cite{Oliynyk2021107} and  \cite{BOOS:2020},  we used  Friedrich's cylinder at infinity conformal transformation to prove the global existence of hyperbolic equations.  While this conformal transformation  works well for wave equations in a  a neighbourhood of space-like infinity,  it  does not work well for wave equations  near future  null  infinity in a neighbourhood of $r=0$.  To address this limitation,  we propose a new mapping that endows the conformal wave equations with the right structure needed to apply the Fuchsian method  near future  null  infinity in a neighbourhood of $r=0$. 

The cylinder at infinity approach  used in \cite{Oliynyk2021107} to compactify  Minnkowski spacetime,  provided insights  suggesting that there exist conformal maps capable of revealing  the structure of the null condition in the conformal spacetime.  Controlling the terms involving the null condition is essential in the proof.  As demonstrated in \cite{Oliynyk2021107},  the associated asymptotic equation to the system of wave equations involves the worst decay terms in a   multiple of the scalar functions  $\bb^{K}_{IJ}$ defined in \cite{Oliynyk2021107}.   The term $\bb^{K}_{IJ}$ is of particular importance since the null condition is satisfied when $\bb^{K}_{IJ}=0$.  From the seminal work of Klainerman,  and Christodoulou,  \cite{Christodoulou:1986,Klainerman1986},  we know that systems  satisfying the null condition admit global solutions.  However,  when $\bb^{K}_{IJ}$ is non-zero,  the null condition is not satisfied, requiring us  to control the decay of the terms involving $\bb^{K}_{IJ}$,  through the bounded weak null condition.

 In \cite{Oliynyk2021107},  the terms $\bb^{K}_{IJ}=0$ can be interpreted as the necessary condition for the null condition to hold in the non-physical bounded manifold.    When the \emph{null condition} is not met,  the function $\bb^{K}_{IJ}$ provide insight into the terms  with the worst decay over time.  Therefore,  identifying the terms $\bb^{K}_{IJ}$ in the non-physical manifold is crucial as their  identification is closely tied to the geometry of the non-physical space under consideration.  At the same time,  this identification is intrinsically related to the conformal map used to transform the physical spacetime.  
 
For wave equations with  quadratic nonlinearities of the form $\ab^{K\mu \nu}_{IJ}\nablab_\mu \ub^I \nablab_\nu \ub^J, $ where  $\ab^{K\mu \nu}_{IJ}$ is a general second order tensor,  the terms $\bb^{K}_{IJ}$ are identified by a Killing vector associated with the conformal transformation.  Specifically,  one can verify that the first column of the Jacobian \eqref{jace} corresponds to a Killing vector in the conformal spacetime \eqref{M+region}.   This identification immediately highlights the terms with the worst decay,   allowing us to control them provided the associated asymptotic system satisfies the bounded weak null condition. Therefore,   conformal transformations of this type,  are strong candidates for our purposes,  as they inherently highlight the terms with the worst decay over time.  This insight can be particularly useful  for classifying conformal transformations that are  suitable for the Fuchsian method. 

\begin{figure}[ht!]
\centering
     \begin{tikzpicture}
          (0,0) node {first node}
            \node [inner sep=0pt,above right]  (plfg)
                {\includegraphics[width=8cm, height=8cm]{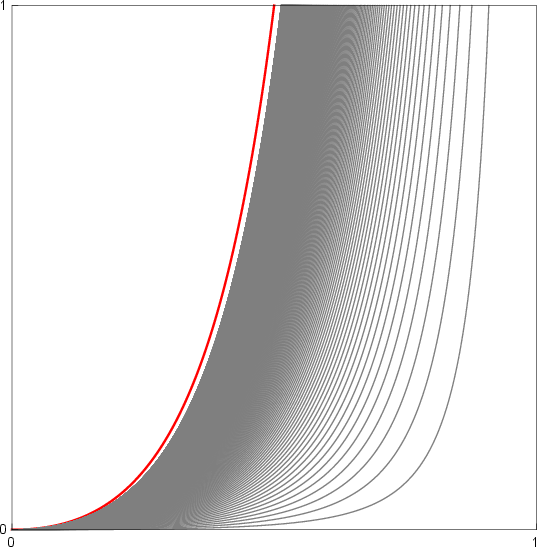}};
         \node (2)[left of=plfg,  xshift=-3.1cm,  yshift=-4cm]{$i^+$};
         \node (3)[below of=plfg,   yshift=-3.3cm]{$\Isc^{+}$};
         \node (4)[right of=plfg,  xshift=3.2cm,    yshift=-3.8cm]{$r$};
         \node (4)[left of=plfg,  xshift=-2.7cm,    yshift=4.3cm]{$t$};
         \node (4)[left of=plfg,  xshift=-0.2cm,    yshift=2.3cm]{$t=\Bigl(\frac{r}{1-r}\Bigr)^2$};
        \end{tikzpicture}
\caption{In this diagram we plot   time-like geodesics of the form $\tb=m\rb, $  from  Minkowski spacetime and represented in the $(t,r,\theta,\phi)$ coordinates,   here  $m \geq 1$.  The red curve  represents the time-like  hyper-surface  $\rb=0$.  In the limit $m \nearrow \infty$ , the time-like curves accumulate near the parabola $t=\Bigl(\frac{r}{1-r}\Bigr)^2$.  Note that all the time like curves end at the point $t=0,\  r=0.$  \label{timep}}
     \begin{tikzpicture}
          (0,0) node {first node}
            \node [inner sep=0pt,above right]  (plfg)
                {\includegraphics[width=8cm, height=8cm]{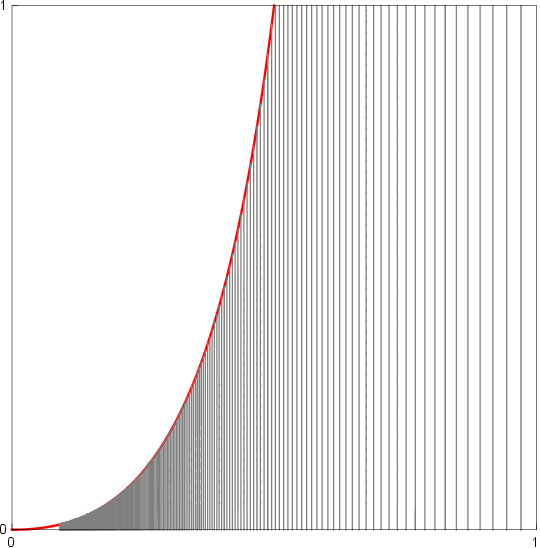}};
         \node (2)[left of=plfg,  xshift=-3.1cm,  yshift=-4cm]{$i^+$};
         \node (3)[below of=plfg,   yshift=-3.3cm]{$\Isc^{+}$};
         \node (4)[right of=plfg,  xshift=3.2cm,    yshift=-3.8cm]{$r$};
         \node (5)[left of=plfg,  xshift=-2.7cm,    yshift=4.3cm]{$t$};
         \node (6)[left of=plfg,  xshift=-0.2cm,    yshift=2.3cm]{$t=\Bigl(\frac{r}{1-r}\Bigr)^2$};
        \end{tikzpicture}
      \caption{In this diagram we plot  null-geodesics $\tb=\rb+b$,  from  Minkowski spacetime and  represented in the $(t,r,\theta,\phi)$ coordinates.   In the limit when $b \nearrow \infty$  the null-geodesics accumulate near $i^+$.   \label{nullp}}
\end{figure}

\begin{figure}
\centering   
     \begin{tikzpicture}
          (0,0) node {first node}
            \node [inner sep=0pt,above right]  (plfg)
                {\includegraphics[width=8cm, height=8cm]{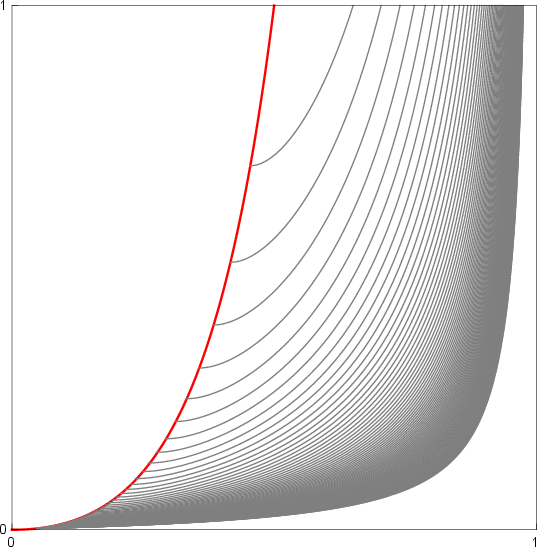}};
         \node (2)[left of=plfg,  xshift=-3.1cm,  yshift=-4cm]{$i^+$};
         \node (3)[below of=plfg,   yshift=-3.3cm]{$\Isc^{+}$};
         \node (4)[right of=plfg,  xshift=3.2cm,    yshift=-3.8cm]{$r$};
         \node (5)[left of=plfg,  xshift=-2.7cm,    yshift=4.3cm]{$t$};
         \node (6)[left of=plfg,  xshift=-0.2cm,    yshift=2.3cm]{$t=\Bigl(\frac{r}{1-r}\Bigr)^2$};
        \end{tikzpicture}
       \caption{In this diagram we plot the family of  space-like geodesics $\tb=k$  from  Minkowski spacetime represented in the $(t,r,\theta,\phi)$ coordinates,  where $k$ is a positive constant and each curve corresponds to a different value of $k$ .   In the limit when $k \nearrow \infty$  the space-like geodesics  accumulate near $\Isc^+$.   \label{spacep}}
     \begin{tikzpicture}
          (0,0) node {first node}
            \node [inner sep=0pt,above right]  (plfg)
                {\includegraphics[width=8cm, height=8cm]{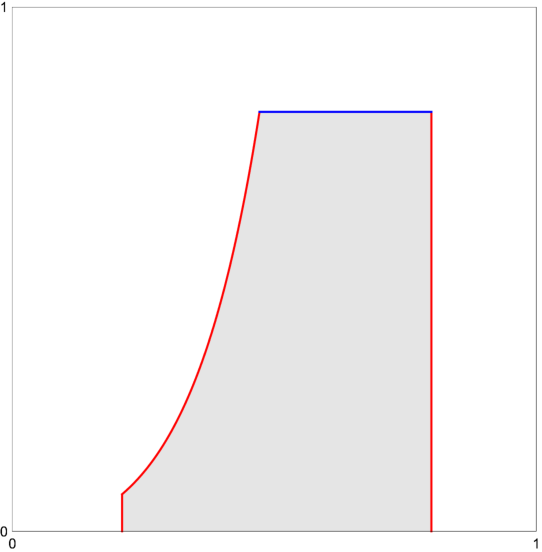}};
         \node (2)[left of=plfg,  xshift=-3.1cm,  yshift=-4cm]{$i^+$};
         \node (3)[below of=plfg,  xshift=2.3cm,  yshift=-3.2cm]{$r_1$};
         \node (4)[right of=plfg,  xshift=3.2cm,    yshift=-3.8cm]{$r$};
         \node (5)[left of=plfg,  xshift=-2.7cm,    yshift=4.3cm]{$t$};
         \node (6)[below of=plfg,  xshift=-2cm,    yshift=-3.2cm]{$r_0$};
         \node (7)[below of=plfg,  xshift=0.8cm  ,yshift=0.2cm]{$M_0$};
         \node (8)[above of=plfg,  xshift=1.1cm  ,yshift=1.7cm]{$\Sigma_0$};
          \node (9)[above of=plfg,  xshift=4.2cm  ,yshift=1.4cm]{$t_0$};
        \end{tikzpicture}
      \caption{The shaded area is the region we are interested in, we prescribe initial data on the  space-like hypersurface $\Sigma_0$. The constants $t_0, r_0, r_1$ will be fixed later in section \ref{extendst}.}
\end{figure}

Using Spherical coordinates $(\xb^{\mu})$ in Minkowski spacetime and the coordinates  $(x^{\mu})$ in the non-physical spacetime,   we define a diffeomorphism  $\psi$ such that
\begin{equation}\label{map}
\psi : \Mb \longrightarrow M : (\xb^i) \longmapsto (x^i) \coloneqq \Bigl(\frac{1}{\tb^2-\rb^2}, \ \frac{1}{1+\tb-\rb},  \theta, \ \phi\Bigr),
\end{equation}
where $\Mb$ is the region \eqref{M+region}, the inverse $\psi^{-1} \ $is given by
\begin{equation*}
\psi^{-1} : M \longrightarrow \Mb : (x^i) \longmapsto (\xb^i) \coloneqq \Bigl(\frac{r^2+t(1-r)^2}{2r(1-r)t}, \ \frac{r^2-t(1-r)^2}{2r(1-r)t},  \theta, \ \phi\Bigr).
\end{equation*}
Note that the  Jacobian  of the map \eqref{map}, is of the form  
\begin{equation}\label{jace}
D\psi(\xb^i)=J^{\mu}_{\alpha}=\frac{\partial x^{\mu}}{\partial \xb_{\alpha}}\bigg\rvert_{\psi^{-1} (x^i)}=\begin{pNiceMatrix}
\frac{-tr^2-(1-r)^2t^2}{r(1-r)}& \frac{tr^2-(1-r)^2t^2}{r(1-r)}&0&0\\
-r^2& r^2&0&0\\
0&0&1&0\\
0&0&0&1
\end{pNiceMatrix},
\end{equation}
the structure of the Jacobian \eqref{jace} ensures that we can factor the terms with lowest decay in time from the components $a^{K00}_{IJ},  \ a^{K01}_{IJ}, a^{K10}_{IJ}, \ a^{K11}_{IJ}$,  (see the expansion \eqref{compa}) as multiples of $b^K_{IJ}$.  This is fundamental to identify the asymptotic equation associated to  \eqref{Mbwave}  since we want the asymptotic equation to \emph{condense} the nonlinearities with the slowest decay.   In  Figures  \ref{timep}, \ref{nullp}, \ref{spacep} we can see the structure of the non-physical spacetime trough the behaviour of its geodesics.  Since we  are interested in a neighbourhood of $i^+$ and $\Isc^+$,  we define  the region $M$ by  
\begin{equation}\notag
\psi(\Mb)=M,
\end{equation}
where $M$ is a non-physical spacetime given by 
\begin{equation}\label{Mregion}
M= \{(t, r) \in (0, \infty)\times (0, 1)    \ |   \ \ \Bigl(\frac{r}{1-r}\Bigr)^2-t>0 \}\times \mathbb{S}^2 . 
\end{equation}
We prescribe  initial data  on the space-like  hyper surface  $\Sigmab$ defined by
\begin{equation}\label{indatS0}
\Sigmab=\bigl\{ (\tb,\rb) \in (0,\infty) \times  (0,\infty)   \  |  \    \tb^2-\rb^2=\frac{1}{t_0},  \ t_0 \in \mathbb{R}^+ \bigr\}  \times \mathbb{S}^2,
\end{equation}
which gets mapped to the non-physical space by
\begin{equation}\label{surfaceSm}
\quad \psi(\Sigmab)= \Sigma. 
 \end{equation} 
We  refine our region of interest  by defining $M_0$ such that 
\begin{equation}\label{M0}
M_{r_0}=\Bigl\{(t, r) \in (0,t_0 )\times \bigl(r_0,  r_1 \bigr)   \Bigl| \quad t < \Bigl(\frac{r}{1-r}\Bigr)^2 \Biggr\}\times \mathbb{S}^2,
\end{equation}
and we restrict the space-like hyper surface \eqref{surfaceSm}  to
\begin{equation}\label{sigmar0}
\Sigma_{0}=\Biggl\{(t, r) \in t_0 \times \Biggl(\frac{t_0^{\frac{1}{2}}}{1+t_0^{\frac{1}{2}}}, r_1\Biggr)  \Biggr\}\times \mathbb{S}^2.
\end{equation}
\newpage
\pagebreak
\newpage
\subsection{Expansion formulas for the tensor components. }
We start with the coordinate transformation
$(\xh^\mu)=(\tb,\rb\cos\phib\sin\thetab,\rb\sin \phib\sin \thetab,\rb\cos \thetab)$,  where $(\hat{x}^{\mu}), \ (\xb^{\mu})$ represent Cartesian and Spherical coordinates respectively,  the Jacobian of this transformation is given by 
\begin{equation} \label{Jbdef}
\Jb^\alpha_\mu =\begin{pmatrix}
 1 & 0 & 0 & 0 \\
 0 & \sin \thetab \cos \phib & \sin \thetab \sin \phib & \cos \thetab \\
 0 & \frac{\cos \thetab \cos \phib}{\rb} & \frac{\cos \thetab \sin \phib}{\rb} & -\frac{\sin \theta}{\rb} \\
 0 & -\frac{\csc \thetab \sin \phib}{\rb} & \frac{\csc \thetab \cos \phib}{\rb} & 0
\end{pmatrix}.
 \end{equation}
Using \eqref{Jbdef}  and the tensorial transformation law
\begin{equation} \label{abcomponents}
\ab_{IJ}^{K\alpha\beta} = \Jb^\alpha_\mu \ah_{IJ}^{K\mu\nu} \Jb^\beta_\nu,
\end{equation}
it is not difficult to verify that  the components of $\ab_{IJ}^{K\alpha\beta}$ can be expanded in powers of $\rb$ as
\begin{equation} \label{abexp}
\ab_{IJ}^{K\alpha\beta} = \cb_{IJ}^{K\alpha\beta}+\frac{1}{\rb}\db_{IJ}^{K\alpha\beta}+ \frac{1}{\rb^2}\eb_{IJ}^{K\alpha\beta},
\end{equation}
where the  coefficients  of \eqref{abexp} are smooth on $\mathbb{S}^2,$ and are classified as follows,  depending on their indices $(\alpha, \beta)$ (see Appendix \ref{indexing} for our indexing conventions): 
\begin{enumerate}
\item[(a)] smooth functions:
$\cb_{IJ}^{K \pc \qc}$,  $\db_{IJ}^{K \pc \qc}$,  $\eb_{IJ}^{K \pc \qc }$,
\item[(b)] smooth vector fields: 
$\cb_{IJ}^{K \qc \Lambda}, \cb_{IJ}^{K \Lambda\qc}$,  $\db_{IJ}^{K \qc \Lambda}$,$\db_{IJ}^{K \Lambda\qc}$,   $\eb_{IJ}^{K \qc \Lambda}$, $\eb_{IJ}^{K \Lambda\qc}$,   
\item[(c)] and smooth (2,0)-tensor fields:  $\cb_{IJ}^{K  \Lambda \Sigma}$,  $\db_{IJ}^{K  \Lambda \Sigma}$, $\eb_{IJ}^{K \Lambda \Sigma }$. 
\end{enumerate}
Explicit formulae for  the components $\cb_{IJ}^{K\alpha\beta}$ can be consulted in \cite{Oliynyk2021107}.  They also can be calculated using \eqref{Jbdef}, and \eqref{abcomponents}.  Note from the definition  \eqref{nullbk},  and the expansion \eqref{abexp}
\begin{equation}\notag
\bb^K_{IJ}=\ab^{K00}_{IJ}-\ab^{K01}_{IJ}-\ab^{K10}_{IJ}+\ab^{K11}_{IJ}=\cb^{K00}_{IJ}-\cb^{K01}_{IJ}-\cb^{K10}_{IJ}+\cb^{K11}_{IJ}.
\end{equation} 
Then we use the tensor transformation
\begin{equation}\label{tentr}
a^{K\mu\nu}_{IJ}=J^{\mu}_{\alpha} \ab^{K\alpha\beta}_{IJ}J^{\nu}_{\beta},
\end{equation}
to write   the tensor $a^K_{IJ}$ in terms of the components $\ab^{K\mu\nu}_{IJ}$,  and the Jacobian $J^{\mu}_{\alpha}$ defined by \eqref{jace}. It is interesting to note  that the components  
\begin{equation}\notag
\{a^{K00}_{IJ}, a^{K01}_{IJ}, a^{K10}_{IJ}, a^{K11}_{IJ}\}
\end{equation}
can be expanded in powers of $t$, such that the lowest order in $t$ has $\bb^K_{IJ}$ as a coefficient. This is due to the particular form of the Jacobian \eqref{jace}.  A  straightforward calculation using \eqref{jace},  and \eqref{tentr}, shows that the components of $a^{K}_{IJ}$ are given by   
\begin{equation}\label{compa}
\begin{split}
a^{K00}_{IJ}=&t^2 \Bigl(\frac{r}{1-r}\Bigr)^2\bb^{K}_{IJ}+2t^3(\cb^{K00}_{IJ}-\cb^{K11}_{IJ})  + \frac{t^4(1-r)^2}{r^2}\bigl(\cb_{IJ}^{K00}+\cb_{IJ}^{K01}+\cb_{IJ}^{K10}+\cb_{IJ}^{K11}\bigr),\\
a^{K01}_{IJ}=&\frac{tr^3}{1-r} \bb^K_{IJ}+ t^2r(1-r)\bigl(\cb_{IJ}^{K00}-\cb_{IJ}^{K01}+\cb_{IJ}^{K10}-\cb_{IJ}^{K11}\bigr),\\
a^{K10}_{IJ}=&\frac{tr^3}{1-r} \bb^K_{IJ}+ t^2r(1-r)\bigl(\cb_{IJ}^{K00}+\cb_{IJ}^{K01}-\cb_{IJ}^{K10}-\cb_{IJ}^{K11}\bigr),\\
a^{K11}_{IJ}=&r^4\bb^{K}_{IJ},\\
a^{K0\Lambda}_{IJ}=&-\frac{2t^2\bigl(r^2+t(1-r)^2\bigr)}{r^2-(1-r)^2t}\db^{K0\Lambda}_{IJ}+2t^2\db^{K1\Lambda}_{IJ},\\
a^{K\Sigma0}_{IJ}=&-\frac{2t^2\bigl(r^2+t(1-r)^2\bigl)}{r^2-(1-r)t}\db^{K\Sigma0}_{IJ}+2t^2\db^{K\Sigma1}_{IJ},\\
a^{K1\Lambda}_{IJ}=&\frac{2r^3(1-r)t}{r^2-(1-r)^2t}\Bigl(\db^{K1\Lambda}_{IJ}-\db^{K0\Lambda}_{IJ}\Bigr),\\
a^{K\Sigma 1}_{IJ}=&\frac{2r^3(1-r)t}{r^2-(1-r)^2t}\Bigl(\db^{K\Sigma 1}_{IJ}-\db^{K\Sigma 0}_{IJ}\Bigr),\\
a^{K\Lambda\Sigma}_{IJ}=&\frac{4r^2(1-r)^2 t^2}{(r^2-(1-r)^2t)^2}\eb^{K\Lambda\Sigma}_{IJ}.
\end{split}
\end{equation}
\section{\bf The conformal wave equation}
Considering the map $\psi$ given by \eqref{map}, we  push-forward the metric \eqref{gbdef}, \eqref{gsldef},  from \eqref{M+region} to \eqref{Mregion}
 \begin{equation}\notag
\gt=\psi_{*}\gb, 
 \end{equation}
the two metrics $\gt$, and $g$  are conformally equivalent and  they satisfy 
\begin{equation}\notag
\gt=\Omega^2g,
\end{equation}
where
\begin{equation}\label{cfact}
\Omega=\frac{r^2-t(1-r)^2}{2r(1-r)t}
\end{equation}
is the conformal factor.   Using \eqref{jace}, \eqref{tentr} we see that  the components of the metric $g$ are given by
\begin{equation}\label{met}
g= -\frac{2r(1-r)}{(r^2-(1-r)^2t)^2}(dt\otimes dr+dr\otimes dt) +\frac{4t}{(r^2-(1-r)^2t)^2}dr\otimes dr +\gsl,
\end{equation}
and  its inverse $g^{-1}$ is given by
\begin{equation}\label{invmet}
g^{-1}= -\frac{(r^2-(1-r)^2t)^2t}{(1-r)^2r^2}\Bigl(\frac{\partial }{\partial t}\otimes \frac{\partial }{\partial t}\Bigr) -\frac{(r^2-(1-r)^2t)^2}{2(1-r)r}\Bigl(\frac{\partial }{\partial t}\otimes \frac{\partial }{\partial r}+\frac{\partial }{\partial r}\otimes \frac{\partial }{\partial t}\Bigr) +\gsl^{-1}.
\end{equation}
We push-forward the wave equation  \eqref{Mbwave} using the map \eqref{map} to obtain
\begin{equation}\label{pform}
\gt^{\mu\nu} \tilde{\nabla}_{\alpha}\tilde{\nabla}_{\beta}\ut^K=\at^{K\mu \nu}_{IJ}\tilde{\nabla}_{\mu}\ut^I \tilde{\nabla}_{\nu}\ut^J, 
\end{equation}
where
\begin{equation}\notag
\ut^K=\psi_{*}\ub^K,  \quad \text{and}   \quad   a^{K\mu\nu}_{IJ}=\psi_* (\ab^{K}_{IJ})^{\mu\nu}.
\end{equation}
The wave equation  \eqref{pform} is well  defined in the region $M_0$ which is given by the  equation \eqref{M0}.  From appendix \eqref{ctrans} and equations \eqref{wavetransA} to \eqref{Nlint} we see that  the wave equation \eqref{pform} is equivalent to  
\begin{equation}\label{weq1}
g^{\alpha\beta} \nabla_{\alpha}\nabla_{\beta}u^K - \frac{n-2}{4(n-1)}Ru^K=f^K,
\end{equation}
where $\nabla$ is  the Levi-Civita connection of $g$ and 
\begin{equation}\label{uom}
u^K=\Omega\ub^K= \frac{r^2-t(1-r)^2}{2r(1-r)t}\ub^K, \quad \quad \text{and} \quad \quad  R=0 .
\end{equation}
Using formulas \eqref{Awet3}, \eqref{Nlint} the source term $f^K$ can be expanded as follows
\begin{equation}\label{fta}
f^K=a^{K\mu \nu}_{IJ} \Bigl(\Omega \nabla_{\mu}u^I \nabla_{\nu}u^J- \nabla_{\mu}u^I\nabla_{\nu}(\Omega) u^J -\nabla_{\nu}u^J\nabla_{\mu}(\Omega)  u^I  +\Omega^{-1}(\nabla_{\mu}\Omega \nabla_{\nu}\Omega) u^I u^J  \Bigr).
\end{equation}
\subsection{\bf The wave equation} 
Let us write explicitly  the wave equation \eqref{weq1} in the coordinates $(x^{\mu})=(t,r,\theta,\phi)$ and by using equations  \eqref{cfact},  \eqref{met},  \eqref{invmet},  \eqref{uom}, we obtain 
\begin{equation}\label{eweq}
t \del{t}\bigl(t\del{t}u^K\bigr)+r(1-r)t\del{r}\del{t}u^K- \frac{r^2(1-r)^2t}{(r^2-(1-r)^2t)^2}\gsl^{\Lambda \Sigma}\nabla_{\Lambda}\nabla_{\Sigma}u^K =- \frac{r^2(1-r)^2t}{(r^2-(1-r)^2t)^2}f^K,
\end{equation}
introducing the change of variables
\begin{equation}\label{chv1}
t\del{t}u^K= (1-r)U_0^K, \quad r\del{r}u^K=t^{-\frac{1}{2}}U_1^K,  \quad \del{\Lambda}u^K= t^{-\frac{1}{2}}(1-r)U_{\Lambda}^K, \quad u^K= t^{-\frac{1}{2}}(1-r)U_4^K, 
\end{equation}
we can transform the wave  equation \eqref{eweq} into the system 
\begin{equation}\label{syst1}
\begin{split}
&\del{t}U_0^K+\frac{(1-r)r}{t}\del{r}U_0^K - \frac{r^2(1-r)^2}{t^{\frac{1}{2}}(r^2-(1-r)^2t)^2}\gsl^{\Lambda \Sigma}\nabla_{\Lambda} U_{\Sigma}^K=\frac{r}{t}U_0^K- \frac{r^2(1-r)}{(r^2-(1-r)^2t)^2}f^K,\\
&\del{t}U^K_1 - \frac{(1-r)r\del{r}U^K_0}{t^{ \frac{1}{2}}}=\frac{U_1^K}{2t}-\frac{r}{t^{\frac{1}{2}}}U^K_0,\\
&\del{t}U^K_{\Lambda} - \frac{\del{\Lambda}U^K_0}{t^{\frac{1}{2}}}=\frac{1}{2t}U^K_{\Lambda}-\frac{r}{t^{\frac{1}{2}}}U^K_0,\\
&\del{t}U^K_4 =\frac{1}{t^{\frac{1}{2}}}U^K_0+\frac{1}{2t}U^K_4.
\end{split}
\end{equation}  
Using   \eqref{cfact},  \eqref{invmet},  \eqref{fta},  \eqref{chv1} we can expand the first term of the right hand side of   \eqref{eweq} 
\begin{equation} \label{bter}
\begin{split}
-&\frac{r^2(1-r)}{(r^2-(1-r)^2t)^2} \Omega a^{K\mu \nu}_{IJ}  \nabla_{\mu}u^I \nabla_{\nu}u^J=-\frac{r}{2t(r^2-(1-r)^2t)}a^{K\mu \nu}_{IJ}  \nabla_{\mu}u^I \nabla_{\nu}u^J   \\
    = &-\frac{r}{2t(r^2-(1-r)^2t)}\Biggl(\frac{a^{K00}_{IJ}}{t^2}(1-r)^2U_0^IU_0^J +\frac{a^{K01}_{IJ}(1-r)}{rt^{\frac{3}{2}}}U_0^IU_1^J +\frac{a^{K10}_{IJ}(1-r)}{rt^{\frac{3}{2}}}U_1^IU_0^J  +\frac{a^{K11}_{IJ}}{r^2t}U_1^IU_1^J  +\\
    &\frac{a^{K0\Lambda}_{IJ}(1-r)^2}{t^{\frac{3}{2}}}U_0^I U_{\Lambda}^J +\frac{a^{K\Sigma0}_{IJ}(1-r)^2}{t^{\frac{3}{2}}}U_{\Sigma}^I U_0^J  + \frac{a^{K1\Lambda}_{IJ}(1-r)}{rt}U_1^I U_{\Lambda}^J  +\\
    &  \frac{a^{K\Sigma 1}_{IJ}(1-r)}{rt}U_{\Sigma}^I U_1^J      + \frac{a^{K\Sigma\Lambda}_{IJ}(1-r)^2}{t}U_{\Sigma}^I U_{\Lambda}^J \Biggr).\\
    \end{split}
\end{equation}
Similarly,  the  second term  of the right hand side of  \eqref{eweq} can be written as
\begin{equation}\notag
    \begin{split}
-&\frac{r^2(1-r)}{(r^2-(1-r)^2t)^2}    a^{K\mu\nu}_{IJ}\Bigl(-\nabla_{\mu}\Omega \nabla_{\nu}u^J u^I -\nabla_{\nu}\Omega \nabla_{\mu}u^I u^J\Bigr)=\\
          &\frac{r^2(1-r)^3}{(r^2-(1-r)^2t)^2}\Biggl[\frac{U_{0}^IU_4^J}{t^{\frac{3}{2}}}\Bigl(-a^{K00}_{IJ}\frac{r}{2(1-r)t^2}  + a^{K01}_{IJ}\frac{r^2+(1-r)^2t}{2r^2(1-r)^2t} \Bigr) +\\
          &\frac{U_{1}^IU_4^J}{r(1-r)t}\Bigl(-a^{K10}_{IJ}\frac{r}{2(1-r)t^2}  + a^{K11}_{IJ}\frac{r^2+(1-r)^2t}{2r^2(1-r)^2t} \Bigr)+\frac{U_{\Lambda}^IU_4^J}{t}\Bigl(-a^{K\Lambda0}_{IJ}\frac{r}{2(1-r)t^2} + a^{K\Lambda1}_{IJ}\frac{r^2+(1-r)^2t}{2r^2(1-r)^2t}  \Bigr)+\\
&\frac{U_{0}^JU_4^I}{t^{\frac{3}{2}}}\Bigl(-a^{K00}_{IJ}\frac{r}{2(1-r)t^2}  + a^{K10}_{IJ}\frac{r^2+(1-r^2)t}{2r^2(1-r)^2t} \Bigr) +\frac{U_{1}^JU_4^I}{r(1-r)t}\Bigl(-a^{K01}_{IJ}\frac{r}{2(1-r)t^2}  + a^{K11}_{IJ}\frac{r^2+(1-r)^2t}{2r^2(1-r)^2t} \Bigr)+\\&\frac{U_{\Sigma}^J U_4^I}{t}\Bigl(-a^{K0\Sigma}_{IJ}\frac{r}{2(1-r)t^2}  + a^{K1\Sigma}_{IJ}\frac{r^2+(1-r)^2t}{2r^2(1-r)^2t}  \Bigr)\Biggr];\\       
    \end{split}
\end{equation}
which after substituting the components of $a^K_{IJ}$ by the expansion \eqref{compa}, simplifies to
\begin{equation}\label{sectr}
    \begin{split}
&\frac{r^2(1-r)^2}{(r^2-(1-r)^2t)^2}\Biggl[\frac{U_{0}^IU_4^J}{rt^{\frac{1}{2}}}\Bigl((-\cb^{K01}_{IJ}+\cb^{K11}_{IJ})r^2-(\cb^{K01}_{IJ}+\cb^{K11}_{IJ})(1-r)^2t \Bigr) +\\
          &\frac{U_{0}^JU_4^I}{rt^{\frac{1}{2}}}\Bigl((-\cb^{K10}_{IJ}+\cb^{K11}_{IJ})r^2-(\cb^{K10}_{IJ}+\cb^{K11}_{IJ})(1-r)^2t \Bigr)+\frac{U_{1}^IU_4^Jr}{t}\bigl(\cb^{K11}_{IJ}-\cb^{K01}_{IJ} \bigr)+\frac{U_{1}^JU_4^Ir}{t}\bigl(\cb^{K11}_{IJ}-\cb^{K10}_{IJ} \bigr)+\\
& U_{\Lambda}^IU_4^J\Bigl(\frac{2r(1-r)^2}{r^2-(1-r)^2t}\db^{K\Lambda1}_{IJ} \Bigr)+U_{\Sigma}^J U_4^I\Bigl(\frac{2r(1-r)^2}{r^2-(1-r)^2t}\db^{K1\Sigma}_{IJ} \Bigr)\Biggr];       
    \end{split}
\end{equation}
and the last term in the right hand side of \eqref{eweq} can be written  as
\begin{equation}\notag
\begin{split}
=-&\frac{2r^3(1-r)^4t}{(r^2-(1-r)^2t)^3}\Biggl(a^{K00}_{IJ}\frac{r^2}{4(1-r)^2t^4} -(a^{K01}_{IJ}+a^{K10}_{IJ})\frac{r^2+(1-r)^2t}{4r(1-r)^3t^3}+a^{K11}_{IJ}\frac{(r^2+(1-r)^2t)^2}{4r^4(1-r)^4t^2}  \Biggr)\frac{U^I_4 U^J_4}{t},\\
\end{split}
\end{equation}
which simplifies to
\begin{equation}\label{thirdfe}
\frac{2r^3(1-r)^4}{(r^2-(1-r)^2t)^3}\cb^{K11}_{IJ}U^I_4U^J_4.
\end{equation}
Substituting into the equations \eqref{bter},  \eqref{sectr},  \eqref{thirdfe},  the  components of $a^K_{IJ}$ given in equation \eqref{compa} and simplifying similar terms it is not difficult to verify that the non-linear  terms from \eqref{syst1} can be expanded as follows
\begin{equation}\notag
\begin{split}
- &\frac{r^2(1-r)^2}{(r^2-(1-r)^2t)^2}f^K=\\
& - \frac{r^3}{2t(r^2-(1-r)^2t)} b^K_{IJ} \Bigl(U_0^I+\frac{U_1^I}{t^{\frac{1}{2}}}\Bigr)  \Bigl(U_0^J+\frac{U_1^J}{t^{\frac{1}{2}}}\Bigr)  +\frac{1}{t}\biggl[ \Bigl( \fc^{K00}_{IJ} U^I_0 U^J_0 +   \fc^{K01}_{IJ} U^I_0 U^J_1  +  \fc^{K10}_{IJ} U^I_1 U^J_0 +  \fc^{K11}_{IJ}U^I_1 U^J_1 +\\ &\fc^{K0\Lambda}_{IJ}U^I_0 U^J_\Lambda+
 \fc^{K\Lambda0}_{IJ}U^I_{\Lambda} U^J_0 + \fc^{K\Lambda \Sigma}_{IJ} U^I_\Lambda U^J_\Sigma+\fc^{K 0\Sigma}_{IJ} U^I_\Lambda U^J_\Sigma +\gc^{K0}_{IJ}U^I_0 U^J_4 + \gc^{K1}_{IJ}U^I_1 U^J_4 + \gc^{K\Lambda}_{IJ} U^I_\Lambda U^J_4 + \hc^{K}_{IJ}U^I_4 U^J_4 
 \Bigr) \biggr],
\end{split}
\end{equation}
where $\{\fc^{K\pc\qc}_{IJ}(t,r)$, $\gc^{K\pc}_{IJ}(t,r)$, $\hc^{K}_{IJ}(t,r)\}$, $\{\fc^{K\pc\Lambda}_{IJ}(t,r)$, $\gc^{K\Lambda}_{IJ}(t,r)\}$
and $\{f^{K\Sigma\Lambda}_{IJ}(t,r)\}$ are collections of smooth scalar, vector,  (2,0)-tensor fields,  respectively on $\mathbb{S}^2$ that
depend smoothly on $(t,r)\in (0,1) \times (0,1)$.  Now consider the following change of variables,  which  is required to  symmetrize the system \eqref{syst1}
\begin{equation}\label{chvm}
V^K=MU^K,
\end{equation}
where
\begin{equation}\label{chvm2}
M=\begin{pmatrix}
t+1&t^{-\frac{1}{2}}&0&0\\
t^{\frac{1}{2}}&\frac{1+\sqrt{5}}{2}\\
0&0&\frac{r(1-r)\mathtt{p}}{(r^2-(1-r)^2t)}&0\\
0&0&0&1\\
\end{pmatrix}, \quad M^{-1}=\begin{pmatrix}
\mathtt{a}_0 &t^{-\frac{1}{2}}\mathtt{a}_1&0&0\\
t^{\frac{1}{2}}\mathtt{a}_1&-(1+t)\mathtt{a}_1&0&0\\
0&0&\frac{(r^2-(1-r)^2t)}{r(1-r)\mathtt{p}}&0\\
0&0&0&1
\end{pmatrix},
\end{equation}
and $\mathtt{a}_0, \mathtt{a}_1, \mathtt{p}$ are smooth  functions of $t$ with 
\begin{equation}\label{a0a1a2}
\mathtt{a}_0=\frac{3+\sqrt{5}}{2+(3+\sqrt{5})t}, \quad \quad \mathtt{a}_1=-\frac{1+\sqrt{5}}{2+(3+\sqrt{5})t},
\end{equation}
and $\mathtt{p}$ is an auxiliary smooth function that we will determine in the next few lines so that  the system in the variables $V^K$ is  symmetric.  Using \eqref{chvm}, \eqref{chvm2}, \eqref{a0a1a2} we write the  system \eqref{syst1} in the form
\begin{align}
\label{syst2}&\del{t}V_0^K+\mathtt{a}_0(1-r)r\del{r}V_0^K+\frac{\mathtt{a}_1}{t^{\frac{1}{2}}}(1-r)r\del{r}V_1^K - \frac{r(1-r)(1+t)}{t^{\frac{1}{2}}\mathtt{p}(r^2-(1-r)^2t)}\gsl^{\Lambda \Sigma}\nabla_{\Lambda} V_{\Sigma}^K=\\\notag
&\mathtt{a}_0(1+r)V_0^K+\frac{\mathtt{a}_1(1+r)}{t^{\frac{1}{2}} }V^K_1- \frac{r^2(1-r)^2}{(r^2-(1-r)^2t)^2}f^K,\\[0.3cm]
\label{syst3}&\del{t}V_1 +\frac{\mathtt{a_1}}{t^{\frac{1}{2}}} (1-r)r\del{r}V^K_0+\frac{2}{t\Bigl(2+(3+\sqrt{5})t\Bigr)}(1-r)r\del{r}V^K_1 - \frac{r(1-r)}{\mathtt{p}(r^2-(1-r)^2t)}\gsl^{\Lambda \Sigma}\nabla_{\Lambda} V_{\Sigma}^K=\\\notag
&\frac{\mathtt{a}_0r(1-\sqrt{5})}{2t^{\frac{1}{2}}}V^K_0+\Bigl(\frac{1}{2t}+\frac{\mathtt{a}_1r(1-\sqrt{5})}{2t}\Bigr)V_1^K,\\[0.3cm]
\label{syst4}&\del{t}V_{\Lambda} - \frac{r(1-r)\mathtt{a}_0\mathtt{p}}{t^{\frac{1}{2}}(r^2-(1-r)^2t)}\nabla_{\Lambda}V^K_0-\frac{r(1-r)\mathtt{a}_1\mathtt{p}}{t(r^2-(1-r)^2t)}\nabla_{\Lambda}V^K_1=\\\notag
&\frac{\mathtt{a}_0(1-r)r^2\mathtt{p}}{t^{\frac{1}{2}}(r^2-(1-r)^2t)}V^K_0+\frac{\mathtt{a}_1(1-r)r^2\mathtt{p}}{t^{\frac{1}{2}}(r^2-(1-r)^2t)}\label{syst5}V^K_1+\Bigl(\frac{r^2+(1-r)^2t}{2t(r^2-(1-r)^2t)}+\frac{\del{t}\mathtt{p}}{\mathtt{p}}\Bigr)V^K_{\Lambda},\\[0.3cm]
&\del{t}V_4 =\frac{\mathtt{a_0}}{t^{\frac{1}{2}}}V^K_0+\frac{\mathtt{a}_1}{t}V_1+\frac{1}{2t}V^K_4.
\end{align}  
To finalize the symetrization of  system \eqref{syst2},  \eqref{syst3}, \eqref{syst4}, \eqref{syst5}  we use the identity $\nabla_{\Lambda} U^K_1=(1-r)r\del{r}U^K_{\Lambda}$,  note that we can write this identity  in the form 
\begin{equation}\label{idq}
\nabla_{\Lambda} U^K_1=(\mathtt{q}+1)(1-r)r\del{r}U^K_{\Lambda}-\mathtt{q}\nabla_{\Lambda}U^K_1,
\end{equation}
where $\mathtt{q}$ is a function that we will determine from a symmetry condition.  Using the change of variable \eqref{chvm}, \eqref{chvm2}, \eqref{a0a1a2},  we can write identity \eqref{idq} as
\begin{equation}\label{idq2}
\nabla_{\Lambda}V_1^K=\frac{t^{\frac{1}{2}}(1+\mathtt{q})}{1+t}\nabla_{\Lambda}V^K_0-\mathtt{q}\nabla_{\Lambda}V^K_1-\frac{(\mathtt{q}+1)(r^2-(1-r)^2t)}{\mathtt{p}(1+t)\mathtt{a}_1}\del{r}V^K_{\Lambda}-\frac{(\mathtt{q}+1)(r^2+(1-r)^2t)}{\mathtt{p}(1+t)\mathtt{a}_1r(1-r)}V^K_{\Lambda},
\end{equation}
substituting \eqref{idq2} into the third equation of \eqref{syst2} we obtain 
\begin{equation}\notag
\begin{split}
\del{t}V^K_{\Lambda}-&\frac{(1-r)r\mathtt{p}}{t^{\frac{1}{2}}(r^2-(1-r)^2t)}\Bigl(\mathtt{a}_0+\frac{\mathtt{a}_1(1+\mathtt{q})}{(1+t)} \Bigr)\nabla_{\Lambda}V_0^K+\frac{ \mathtt{a}_1\mathtt{p} \mathtt{q}(1-r)r}{t(r^2-(1-r)^2t)}\nabla_{\Lambda}V^K_1+\frac{(\mathtt{q}+1)}{t(1+t)}(1-r)r\del{r}V^K_{\Lambda}=\\
&\frac{\mathtt{a}_0(1-r)r^2\mathtt{p}}{t^{\frac{1}{2}}(r^2-(1-r)^2t)}V^K_0+\frac{\mathtt{a}_1(1-r)r^2\mathtt{p}}{t^{\frac{1}{2}}(r^2-(1-r)^2t)}V^K_1+\Bigl(\frac{(r^2+(1-r)^2t)(3+t+2q)}{2t(r^2-(1-r)^2t)(1+t)}+\frac{\del{t}\mathtt{p}}{\mathtt{p}}\Bigr)V^K_{\Lambda},
\end{split}
\end{equation}
the system \eqref{syst2}- \eqref{syst5}  is symmetric  if we impose the condition 
\begin{equation}\label{systpq}
\mathtt{p}\mathtt{a}_0+\frac{\mathtt{a}_1\mathtt{p}(1+\mathtt{q})}{(1+t)} =\frac{1+t}{\mathtt{p}} \quad  \text{and} \quad  \frac{\mathtt{a}_1\mathtt{p} \mathtt{q}}{t}=-\frac{1}{\mathtt{p}};
\end{equation}
solving the system \eqref{systpq} for $\mathtt{p}, \mathtt{q}$ we obtain
\begin{equation}\label{pq}
\mathtt{p}=\sqrt{1+t(3+t)},  \quad \text{and} \quad \mathtt{q}=\frac{\left(1+\sqrt{5}\right) t}{3+\sqrt{5}+2 t}.
\end{equation}
Then,  we proceed to write the symmetric system  \eqref{syst2},  \eqref{idq2},  \eqref{pq}  in the form
\begin{equation}\label{shweq}
B^0\del{t}V^K+\frac{1}{t}B^1(1-r)r\del{r}V^K +\frac{1}{t^{\frac{1}{2}}}B^{\Lambda}\nabla_{\Lambda}V^K=\frac{1}{t}\mathcal{B}\mathbb{P}V^K+\frac{1}{t^{\frac{1}{2}}}\mathcal{C}V^K+F^K,
\end{equation}
where
\begin{align}
\label{VKvar} V^K&=(V^K_{\mathcal{I}})= \begin{pmatrix}V^K_0\\ V^K_1 \\ V^K_{\Lambda} \\ V^K_4 \end{pmatrix}, \\[0.9cm]
\label{MB0}B^0&=\begin{pmatrix}
1&0&0&0\\
0&1&0&0\\
0&0 &\delta_{\Omega}^{\Sigma} & 0\\
0&0&0&1
\end{pmatrix},\\[0.9cm]
\label{MB1}B^1&=\begin{pmatrix}
t \mathtt{a}_0 & t^{\frac{1}{2}}\mathtt{a}_1 & 0 & 0\\
t^{\frac{1}{2}}\mathtt{a}_1 & \frac{2}{2+(3+\sqrt{5})t}&0&0\\
0&0&\delta^{\Sigma}_{\Omega}\frac{2}{2+(3-\sqrt{5})t}&0\\
0&0&0&0\\
\end{pmatrix},  \\[0.9cm]
\label{MBL}B^{\Lambda}&=\begin{pmatrix}
0&0&-\gsl^{\Lambda\Sigma}\frac{(1-r)r(1+t)}{\mathtt{p}(r^2-(1-r)^2t)}&0\\
0&0&-\gsl^{\Lambda\Sigma}\frac{t^{\frac{1}{2}}(1-r)r}{\mathtt{p}(r^2-(1-r)^2t)}&0\\
-\delta^{\Lambda}_{\Omega}\frac{(1-r)r(1+t)}{\mathtt{p}(r^2-(1-r)^2t)}&-\delta^{\Lambda}_{\Omega}\frac{t^{\frac{1}{2}}(1-r)r}{\mathtt{p}(r^2-(1-r)^2t)}&0&0\\
0&0&0&0
\end{pmatrix},  \\[0.9cm]
\label{MBB}\mathcal{B}&=\begin{pmatrix}
\frac{1}{2}&0&0&0\\
0&\frac{1}{2}+ \frac{\mathtt{a}_1r(1-\sqrt{5})}{2}&0&0\\
0&0&\delta^{\Lambda}_{\Omega}\frac{(r^2+(1-r)^2t)(3+t+2\mathtt{q})}{2(r^2-(1-r)^2t)(1+t)}&0\\
0&\mathtt{a}_1&0&\frac{1}{2}\\
\end{pmatrix},  \\[0.9cm]
\label{MBC}\mathcal{C}&=\begin{pmatrix}
t^{\frac{1}{2}}\mathtt{a}_0(1+r)&\mathtt{a}_1(1+r)&0&0\\
\frac{\mathtt{a}_0r(1-\sqrt{5})}{2}&0&0&0\\
\frac{\mathtt{a}_0(1-r)r^2\mathtt{p}}{r^2-(1-r)^2t}&\frac{\mathtt{a}_1(1-r)r^2\mathtt{p}}{r^2-(1-r)^2t}&\delta^{\Lambda}_{\Omega}\frac{(3+2t)t^{\frac{1}{2}}}{2+2t(3+t)}&0\\
\mathtt{a}_0&0&0&0\\
\end{pmatrix},  \\[0.9cm]
\mathbb{P}&=\begin{pmatrix} \label{Pop}
0&0&0&0\\
0&1&0&0\\
0&0&\delta^{\Lambda}_{\Omega}&0\\
0&0&0&1
\end{pmatrix}, \\[0.9cm]
\notag F^K&=\begin{pmatrix}-\frac{(1+t)(1-r)r^2}{(r^2-(1-r)^2t)^2} f^K \\  -\frac{t^{\frac{1}{2}}(1-r)^2r^2}{(r^2-(1-r)^2t)^2} f^K \\ 0 \\ 0 \\ 0 \end{pmatrix}.
\end{align}
The operator  $\Pbb$ is a covariantly constant, time-independent, symmetric projection operator that commutes
with $B^0$,  and $\Bc$, that is,
\begin{equation} \label{Poppr1}
\Pbb^\perp = \id - \Pbb,  \quad \Pbb^2 = \Pbb, \quad \Pbb^{\tr}=\Pbb, \quad \quad \del{t}\Pbb=0, \quad \del{r}\Pbb =0, \AND \nablasl{\Lambda}\Pbb = 0,
\end{equation} 
and
\begin{equation} \label{Poppr2}
[B^0,\Pbb] = [\Bc,\Pbb]=0,
\end{equation}
where the symmetry is with respect to the inner-product
\begin{equation} \label{hdef}
h(Y,X) = \delta^{\pc\qc}Y_\pc X_\qc + \gsl^{\Sigma\Lambda}Y_\Lambda X_\Sigma + Y_4 X_4.
\end{equation}
Moreover,  note from the definitions \eqref{MB0}, \eqref{MB1}, \eqref{MBL}  that $B^0$,  $B^1$ and $B^\Lambda$ are symmetric with respect to \eqref{hdef}
and that $B^0$ satisfies
\begin{equation} \notag
h(Y,Y)= h(Y,B^0Y),
\end{equation} 
for all $Y=(Y_\Ic)$ and $0<t\leq t_0$,  which implies that the system \eqref{shweq} is symmetric hyperbolic.  Now  we introduce a change of radial coordinate via
\begin{equation}\label{rtr}
r=\rho^m, \quad \quad m\in \mathbb{Z}_{m \geq 1},
\end{equation}
and note that the  transformation \eqref{rtr} leads to
\begin{equation}\label{rhopr}
r\del{r} = \frac{\rho}{m}\del{\rho},
\end{equation}
then,  using \eqref{rtr}, \eqref{rhopr} we can express \eqref{shweq} as
\begin{equation}\label{shweq2}
B^0\del{t}V^K+\frac{1}{t}B^1\frac{(1-\rho^m)\rho}{m}\del{\rho}V^K +\frac{1}{t^{\frac{1}{2}}}B^{\Lambda}\nabla_{\Lambda}V^K=\frac{1}{t}\mathcal{B}\mathbb{P}V^K+\frac{1}{t^{\frac{1}{2}}}\mathcal{C}V^K+F^K,
\end{equation}
where now any $r$ appearing in $B^{\Lambda},\mathcal{B},  F^K$ is replaced using \eqref{rtr}.  Notice that without loosing generality we can choose new constants $r_0^{\frac{1}{m}}=\rho_0, r_1^{\frac{1}{m}}=\rho_1$ which define a new  spacetime region \eqref{M0} expressed in terms of the radial coordinate $\rho$ as
\begin{equation}\label{M02}
M_{\rho_0}=\Biggl\{(t, \rho) \in \bigl(0, t_0\bigr)\times (\rho_0, \rho_{1})  \  \Biggl| \  
t < \frac{\rho^{2m}}{(1+\rho^{m})^{2}},  \ \rho_0,\rho_1 \in (0,1) \Biggr\}\times \mathbb{S}^2,
\end{equation}
and we redefine  the space-like hyper surface \eqref{sigmar0}  were we prescribe initial data as
\begin{equation}\label{srho0}
\Sigma= \Biggl\{(t,\rho) \in t_0 \times \Biggl(\frac{t_0^{\frac{1}{2m}}}{(1+t_0^{\frac{1}{2}})^{\frac{1}{m}}}, \rho_1\Biggr) \Biggr\} \times \mathbb{S}.
\end{equation}
\subsection{The extended system }\label{extendst}
Next, we let $\alpha>0$ and $\chi(\rho)$ denote a constant,  and smooth cut-off function that satisfies
\begin{equation*}
\chi \geq 0, \quad \chi|_{[\rho_0,\rho_1]} = 1,  \AND \supp(\chi) \subset (\rho_0-\alpha,\rho_1+\alpha),
\end{equation*}
we then consider an extended version of \eqref{shweq2} given by
\begin{equation} \label{MwaveE}
B^0 \del{t} V^K +  \frac{1}{t}\frac{\chi(1-\rho^m)\rho }{m} B^1\del{\rho} V^K + \frac{\chi}{t^{\frac{1}{2}}}B^\Lambda \nablasl{\Lambda} V^K = \frac{1}{t}\tilde{\Bc} \mathbb{P}V^K +\frac{1}{t^{\frac{1}{2}}}\tilde{\mathcal{C}}V^K+\mathcal{F}^K  
\end{equation} 
where
\begin{gather} 
\mathcal{F}^K = \frac{1}{t}Q^K \ev_0 + \Gc^K, \label{FcKdef} \\
Q^K = -\frac{\rho^{3m}}{2\bigl(\rho^{2m}-(1-\rho^m)^2t\bigr)}b^{K}_{IJ}\chi(\rho) V^I_0 V_0^J, \label{QKdef} \\
\Gc^K = 
\Gc_0+\frac{1}{t^{\frac{1}{2}}}\Gc_1 +
 \frac{1}{t}\Gc_2, \label{Gcdef}\\
\Gc_0^K =   G_0^K(t^{\frac{1}{2}},t,\chi(\rho)\rho^m,V,V), \label{GcK0def} \\
\Gc_1^K =  G_1^K(t^{\frac{1}{2}},t,\chi(\rho)\rho^m,V,\Pbb V), \label{GcK1def}\\
\Gc_2^K =  G_2^K(t^{\frac{1}{2}},t,\chi(\rho)\rho^m,\Pbb V,\Pbb V), \label{GcK2def}
\end{gather}
where the $G^K_a(\tau,t,\rho, V,\hat{V})$ are smooth bilinear maps,  that is
\begin{equation*}
 G^K(\tau,t,\rh, X,Y)=G^{K\qc\pc}_{IJ}(\tau,t,\rho)V^I_{\qc}\hat{V}^J_{\pc}+G^{K\qc\Lambda}_{IJ}(\tau,t,\rho)V^I_{\qc}\hat{V}^J_{\Lambda}+G^{K\Lambda\Sigma}_{IJ}(\tau,t,\rho)V^I_{\Lambda}\hat{V}^J_{\Sigma}
 \end{equation*} 
 corresponding to  smooth scalar,  vector, and $(0,2)$ tensor fields respectively,  that depend smoothly on $(\tau,t,r)$,  and
 \begin{equation}\label{PbbGc2=0}
 \mathbb{P}G^K_2=0.
 \end{equation}
 The maps $\tilde{\Bc}, \tilde{\Cc}$ are defined by
\begin{align}
\label{BMEx2}\tilde{\Bc}=&\Bc_{*}+ \chi (\Bc-\Bc_{*}),\\
\label{BMEx3}\tilde{\Cc}=&\Cc_{*}+ \chi (\Cc-\Cc_{*}),
\end{align}
where we are using the notation
\begin{equation}\notag
(\cdot)_{*}=(\cdot)\bigl|_{\rho_1=1}.
\end{equation}
Note that the system  \eqref{MwaveE} is well-defined on the extended spacetime region
\begin{equation}\label{extst}
(0, t_0) \times\mathcal{S},  
\end{equation} 
with
\begin{equation}\label{extst2}
\mathcal{S}= T^{1}_{\alpha}\times \mathbb{S}^2,
\end{equation}
 and $T^{1}_{\alpha}$ is the 1-dimensional torus obtained from identifying the end points of the interval $[\rho_0-2\alpha,\rho_1+2\alpha]$.  We determine the  value of $t_0$ in the calculations below.   By construction, \eqref{MwaveE} agrees with \eqref{shweq2} when restricted to \eqref{M02}.  Evaluating   \eqref{MBB}, \eqref{MBC} at $\rho_1=1$ we obtain
\begin{align}
\label{BCS1} \Bc_{*}&=\begin{pmatrix}
\frac{1}{2}&0&0&0\\
0&\frac{1}{2}+ \frac{\mathtt{a}_1(1-\sqrt{5})}{2}&0&0\\
0&0&\delta^{\Lambda}_{\Omega}\frac{3+t+2\mathtt{q}}{2(1+t)}&0\\
0&\mathtt{a}_1&0&\frac{1}{2}\\
\end{pmatrix},\\
\label{BCS2}\mathcal{C}_{*}&=\begin{pmatrix}
2t^{\frac{1}{2}}\mathtt{a}_0&2\mathtt{a}_1&0&0\\
\frac{\mathtt{a}_0(1-\sqrt{5})}{2}&0&0&0\\
0&0&\delta^{\Lambda}_{\Omega}\frac{(3+2t)t^{\frac{1}{2}}}{2+2t(3+t)}&0\\
\mathtt{a}_0&0&0&0\\
\end{pmatrix}.
\end{align}
To proceed,  note that  $B^0=\mathbbm{1}$ is positive definite and  
\begin{equation}\label{InpV}
h(V,B^0V)= h(V,V),
\end{equation}
this makes  clear  that the system \eqref{MwaveE} is symmetric hyperbolic on $(0, t_0)\times T^1_{\alpha}\times \mathbb{S}^2$.  Now we verify that the operators  $B^{\mu}, B^{\nu},$  $\mathcal{B}, \tilde{\mathcal{B}},  \mathcal{C},  \tilde{\mathcal{C}}, $ are bounded in \eqref{extst}.   First we set  $m \in \mathbb{N}$ and $0<\eta <1 $,  then using Taylor's theorem,  equations  \eqref{MB0}-\eqref{MBC}  and  \eqref{BCS1},  \eqref{BCS2},  it can be verified that there exist a constant such that
\begin{equation*}
C(m, \eta)>0 
\end{equation*} 
and
\begin{align}
\label{BC1}|\mathcal{B}-\mathcal{B}_*|&\leq C|\rho|,\\
\label{BC2}|\del{\rho}\mathcal{B}|&\leq C, \\
\label{BC3}|\mathcal{C}-\mathcal{C}_*|&\leq C|\rho|,\\
\label{BC4}|\del{\rho}\mathcal{C}|&\leq C, \\
\notag|B^{\Lambda}|&\leq C|\rho|,\\
\notag|\del{\rho}B^{\Lambda}|&\leq C, 
\end{align}
for all $(t, \rho, x^{\Lambda}) \in (0,t_0)\times (\rho_0-\alpha, \eta)\times \mathbb{S}^2$.  Therefore, using \eqref{BC1},  there exist a constant $\sigma>0$ that  bounds 
\begin{equation}\label{BC7}
|\tilde{\Bc}-\Bc_*| = |\chi| |\Bc-\Bc_*| \leq |\Bc-\Bc_*|< C |\rho| < \sigma.
\end{equation}
Applying a similar reasoning and using equation \eqref{BC3},  we see that 
\begin{equation}\label{BC9}
|\tilde{\Cc}-C_*|=|\chi||\Cc-\Cc_*| \leq |\Cc-\Cc_*|<C|\rho|<\sigma,
\end{equation}
for all $(t, \rho, x^{\Lambda}) \in (0,t_0)\times (\rho_0-\alpha,\rho_1+\alpha)\times \mathbb{S}^2$,  where 
\begin{equation}\notag
0<\rho_1-\rho_0+2\alpha < \min \left\{\eta, \sigma \right\}.
 \end{equation} 
Using  \eqref{MB0},  \eqref{MB1},  we also define  $\sigma_1>0$ such that
\begin{equation}\label{supB}
\sigma_1= \sup_{t\in(0,  t_0)}\left\{|B^0|,  |B^1|\right\}.
\end{equation}
We then proceed  by nothing that
 \begin{equation}\notag
 \del{\rho}\Bigl(\frac{(1-\rho^m)\rho}{m}B^1\Bigr)= \del{\rho}\chi(\rho)\frac{(1-\rho^m)\rho}{m}B^1 +\chi(\rho)\frac{1-\rho^m(1+m)}{m}B^1,
 \end{equation}
 therefore,  using  definition \eqref{supB} we can write the above equation in the form
 \begin{equation}
 \begin{split}
\label{boundpB1} \left\lvert \del{\rho}\Bigl(\frac{(1-\rho^m)\rho}{m}B^1\Bigr)\right \lvert &\leq \lVert\del{\rho} \chi(\rho)\rVert_{L^{\infty}(\mathbb{R})} \Bigl\lvert \frac{(1-\rho^m)\rho}{m}\Bigr\rvert | B^1| +|\chi(\rho)| \Big\lvert\frac{1-\rho^m(1+m)}{m}\Bigr\lvert |B^1|,\\
\left\lvert \del{\rho}\Bigl(\frac{(1-\rho^m)\rho}{m}B^1\Bigr)\right \lvert&< \lVert\del{\rho} \chi(\rho)\rVert_{L^{\infty}(\mathbb{R})}\frac{\sigma_1}{m} +\frac{1+m}{m}\sigma_1.
\end{split}
 \end{equation}
Similarly,  using \eqref{BC1},  \eqref{BC2},  \eqref{BC7},  we see that
 \begin{equation}\notag
 \del{\rho}\tilde{\mathcal{B}}=\del{\rho}\chi(\rho)\left( \Bc -\Bc_* \right)+\chi(\rho) \del{\rho}\Bc,
 \end{equation}
 so we  conclude that
 \begin{equation}\notag 
 | \del{\rho}\tilde{\mathcal{B}}|< \left(\lVert\del{\rho} \chi(\rho)\rVert_{L^{\infty}(\mathbb{R})}+1\right)\sigma,
 \end{equation}
 applying  similar arguments and using \eqref{BC3},  \eqref{BC4}, \eqref{BC9} it is not difficult to show that 
 \begin{equation}\notag 
 | \del{\rho}\tilde{\mathcal{C}}|< \left(\lVert\del{\rho} \chi(\rho)\rVert_{L^{\infty}(\mathbb{R})}+1\right)\sigma,
 \end{equation}
for all $(t,\rho,x^{\Lambda})\in (0, t_0)\times(0,1)\times\mathbb{S}^2$.  Note from definition \eqref{M02} that the boundary of the region $M_{\rho_0}$  can be decomposed as
\begin{equation*}
\del{}M_{\rho_0} = \Sigma \cup \Sigma^+ \cup \Gamma^{-} \cup \Gamma^{+}\cup \Gamma,
\end{equation*}
where
\begin{align}
\label{BSigma}\Sigma&=\{t_0\}\times \Biggl(\frac{t_0^{\frac{1}{2m}}}{(1+t_0^{\frac{1}{2}})^{\frac{1}{m}}},  1 \Biggr)\times \mathbb{S}^2,\\
\label{BSigma+} \Sigma^+ &= \{0\}\times \bigl(\rho_0,1\bigr)\times \mathbb{S}^2, \\
\Gamma^{-} &= \Biggl[0,\frac{\rho_0^{2m}}{(1-\rho_0^m)^2} \Biggr] \times \{\rho_0\} \times \mathbb{S}^2,\label{n1}\\
\Gamma^{+} &= [0,t_0]\times \{1\}  \times \mathbb{S}^2,\label{n2}\\
\label{BGamma}\Gamma &=\Bigl\{(t,\rho) \in \Bigl[\frac{\rho_0^{2m}}{(1-\rho_0^m)^2}, t_0\Bigr]\times \Bigl[\rho_0,\frac{t_0^{\frac{1}{2m}}}{(1+t_0^{\frac{1}{2}})^{\frac{1}{m}}}\Bigr] \quad \Bigl| \quad  t=\frac{\rho^{2m}}{(1-\rho^m)^2} \Bigr\}\times \mathbb{S}^2,
\end{align} 
where \eqref{BSigma} is the space-like hypersurface where we prescribed initial data and \eqref{BSigma+} corresponds to a section of $\Isc^{+}$.  In the limit when $\rho_0 \searrow 0$ then,  $\Sigma^{+}$ corresponds with $\Isc^{+}$.  Using \eqref{n1}, \eqref{n2}, \eqref{BGamma} we calculate the following co-normals
\begin{equation}\label{nmls}
n^{-}= -d\rho,  \quad  \ n^{+}=d\rho,  \quad  \text{and} \quad n = -dt + \frac{2m\rho^{2m-1}}{(1-\rho^{m})^3} d\rho, 
\end{equation}
which  define outward pointing co-normals to $\Gamma^{-}$,  $\Gamma^{+}$,  and $\Gamma$ respectively.  Furthermore, we have from \eqref{MB0}, \eqref{MB1}, \eqref{MBL} and \eqref{nmls} that
\begin{align}
\label{wspacelike1}
\Bigl(n^{+}_{0} B^0 +n^{+}_1 \frac{\chi(1-\rho^m)\rho}{m} B^1 + n^{+}_\Sigma B^\Sigma\Bigr) \Bigl|_{\Gamma^{+}} &= 0, \\
\label{wspacelike2}  \Bigl(n^{-}_{0} B^0 +n^{-}_1 \frac{\chi(1-\rho^m)\rho}{m} B^1 + n^{-}_\Sigma B^\Sigma\Bigr)\Bigl|_{\Gamma^{-}}  &= - \frac{(1-\rho_0^m)\rho_0}{m} B^1\Bigl|_{\Gamma^{-}},\\
\label{wspacelike3} \Bigl(n_{0} B^0 +n_1 \frac{\chi(1-\rho^m)\rho}{m} B^1 + n_\Sigma B^\Sigma\Bigr)\Bigl|_{\Gamma}  &= -B^0+2tB^1\Bigl|_{\Gamma},
\end{align}
where in deriving \eqref{wspacelike3} we have used  \eqref{BGamma} with $t=\frac{\rho^{2m}}{(1-\rho^m)^2}$ on $\Gamma$.  From inequalities \eqref{wspacelike1}, \eqref{wspacelike2}, \eqref{wspacelike3} 
we deduce that
\begin{equation}\label{iny2}
h\Bigl(Y,\Bigl(n^{+}_{0} B^0 +n^{+}_1 \frac{\chi(1-\rho^m)\rho}{m} B^1 + n^{+}_\Sigma B^\Sigma\Bigr)\Bigl|_{\Gamma^{+}}Y\bigr) = 0\\
\end{equation}
and
\begin{equation}
\begin{split}\label{iny0}
& h\Bigl(Y, \Bigl(n^{-}_{0} B^0 +n^{-}_1 \frac{\chi(1-\rho^m)\rho}{m} B^1 + n^{-}_\Sigma B^\Sigma\Bigr)\Bigl|_{\Gamma^{-}} Y\Bigr)=\\
&-\frac{(1-\rho_0^m)\rho_0}{m}\Bigl(t\mathtt{a}_0Y_0^2 +2t^{\frac{1}{2}}\mathtt{a}_1Y_0Y_1 +\frac{2}{2+(3+\sqrt{5})t}Y_1^2 +\frac{2}{2+(3-\sqrt{5})t}\gsl^{\Lambda\Sigma}Y_{\Lambda}Y_{\Sigma} \Bigr)\leq \\
&-\frac{(1-\rho_0^m)\rho_0}{m}\Biggl(t\Bigl(\mathtt{a}_0+\mathtt{a}_1\beta_0\Bigr)Y_0^2 +\Bigl(\frac{2}{2+(3+\sqrt{5})t}+ \frac{\mathtt{a}_1}{\beta_0}\Bigr)Y_1^2 +\frac{2}{2+(3-\sqrt{5})t}\gsl^{\Lambda\Sigma}Y_{\Lambda}Y_{\Sigma} \Biggr) 
\end{split}
\end{equation}
where the constant $\beta_0$ in \eqref{iny0} comes from  Young's inequality.\footnote{Here we use Young's inequality in the form  $ |ab|\leq\frac{a^2\beta_0}{2}+\frac{b^2}{2\beta_0}$,  thus  $-\frac{t\beta_0Y_0^2}{2}-\frac{Y_1^2}{2\beta_0} \leq t^{\frac{1}{2}}Y_0Y_1 \leq \frac{tY_0^2\beta_0}{2}+\frac{Y_1^2}{2\beta_0}$.  Since $\mathtt{a}_1<0$ we get that $-t\beta_0\mathtt{a}_1Y_0^2-\frac{\mathtt{a}_1Y_1^2}{\beta_0} \geq 2t^{\frac{1}{2}}\mathtt{a}_1Y_0Y_1$ .} Choosing $\beta_0=\frac{1}{2}(1+\sqrt{5})$ implies  $\mathtt{a}_0+\mathtt{a}_1\beta_0=\frac{2}{2+(3+\sqrt{5})t}+ \frac{\mathtt{a}_1}{\beta_0}=0$,  which leads to 
\begin{equation}\label{iny}
\begin{split}
h\Bigl(Y, \Bigl(n^{-}_{0} B^0 +n^{-}_1 \frac{\chi(1-\rho^m)\rho}{m} B^1 + n^{-}_\Sigma B^\Sigma\Bigr)\Bigl|_{\Gamma^{-}} Y\Bigr)\leq-\frac{(1-\rho_0^m)\rho_0}{m}\Biggl(\frac{2}{2+(3-\sqrt{5})t}\gsl^{\Lambda\Sigma}Y_{\Lambda}Y_{\Sigma} \Biggr)
\end{split}\leq0
\end{equation}
for all $(t,\rho)\in \Gamma^{-}$. We proceed in a similar way with inequality \eqref{wspacelike3} to obtain
\begin{equation}\notag
\begin{split}
h\Bigl(Y, \Bigl(n_{0} B^0 +n_1 \frac{\chi(1-\rho^m)\rho}{m} B^1 + n_\Sigma B^\Sigma\Bigr)\Bigl|_{\Gamma} Y\Bigr)=h\Bigl(Y, \Bigl(-B^0+2tB^1\Bigr)\Bigl|_{\Gamma}Y \Bigr)=\\
(-1+2t^2\mathtt{a}_0)Y_0^2+4t^{\frac{3}{2}}\mathtt{a}_1Y_0Y_1 +\Bigl(-1+\frac{4\mathtt{a}_0t}{3+\sqrt{5}}\Bigr)Y_1^2+\Bigl(-1+\frac{4t}{2+(3-\sqrt{5})t}\Bigr)\gsl^{\Lambda \Sigma} Y_{\Lambda}Y_{\Sigma}\leq\\
(-1+2t^2\mathtt{a}_0-2\mathtt{a}_1t^3\beta_1)Y_0^2+\Bigl(-1+\frac{4\mathtt{a}_0t}{3+\sqrt{5}}-\frac{2\mathtt{a}_1}{\beta_1}\Bigr)Y_1^2+\Bigl(-1+\frac{4t}{2+(3-\sqrt{5})t}\Bigr)\gsl^{\Lambda \Sigma} Y_{\Lambda}Y_{\Sigma}
\end{split}
\end{equation}
 where the constant $\beta_1$comes from  Young's inequality.  Choosing $\beta_1 =1+\sqrt{5}$,  implies   $-1+2t^2\mathtt{a}_0-2\mathtt{a}_1t^3\beta_1\leq 0$,  and  $ -1+\frac{4\mathtt{a}_0t}{3+\sqrt{5}}-\frac{2\mathtt{a}_1}{\beta_1} \leq 0$ for all $t \in (0, t_0)$ where we set
 \begin{equation}\label{time0}
 t_0=\frac{1}{8} \Bigl( \sqrt{5} + \sqrt{ 10 \sqrt{5}-2}-3\Bigr) 
  \end{equation} and nothing that $-1+\frac{4t}{2+(3-\sqrt{5})t}\leq0$ for $t \in (0, t_0)$ we conclude that
\begin{equation}\label{iny3}
h\Bigl(Y, \Bigl(n_{0} B^0 +n_1 \frac{\chi(1-\rho^m)\rho}{m} B^1 + n_\Sigma B^\Sigma\Bigr)\Bigl|_{\Gamma} Y\Bigr)=h\Bigl(Y, \Bigl(-B^0+2tB^1\Bigr)Y \Bigr)\leq0
\end{equation}
for $t\in(0,t_0)$.  From equations \eqref{iny2}, \eqref{iny}, \eqref{iny3} and  the definition given by \cite[\S 4.3]{Lax:2006},  the surfaces $\Gamma^{+}$,  $\Gamma^-$  and $\Gamma$, are weakly spacelike  for all $Y=(Y_\Ic)$ and $t\in(0,t_0)$.  Note that    solution of the extended system \eqref{MwaveE} on the extended spacetime \eqref{extst} yields  a solution of the original system  \eqref{shweq2} on the region \eqref{M02}.  The  solution  is uniquely determined by the restriction of the initial data to \eqref{srho0}.  To continue first we must verify a structural condition regarding the  operators $\tilde{\mathcal{B}}, B^0$.  Using \eqref{hdef} and \eqref{BCS1} we have
\begin{equation}\label{h0B}
h\left(V, \mathcal{B}_* V\right)=\delta_{IJ}\left(\frac{1}{2}V_0^I V_0^J + \left(\frac{1}{2}+\mathtt{a}_1\frac{1-\sqrt{5}}{2}\right)V_1^I V_1^J +\left(\frac{3+t+2\mathtt{q}}{2(1+t)}\right)\gsl^{\Lambda \Sigma}V_{\Lambda}^IV_{\Sigma}^J +\mathtt{a}_1Y_4^I Y^J_1+\frac{1}{2}Y_4^I Y^J_4\right),
\end{equation}
using a similar version of Young's inequality from \eqref{iny0},  we can write \eqref{h0B} as follows
\begin{equation}\label{h1B}
\begin{split}
h\left(V, \mathcal{B}_* V\right)\geq\delta_{IJ}\Biggl(\frac{1}{2}V_0^I V_0^J + \left(\frac{1}{2}+\mathtt{a}_1\frac{1-\sqrt{5}}{2}+ \frac{\mathtt{a}_1\beta_2}{2}\right)V_1^I V_1^J +\\
\left(\frac{3+t+2\mathtt{q}}{2(1+t)}\right)\gsl^{\Lambda \Sigma}V_{\Lambda}^IV_{\Sigma}^J +\left(\frac{1}{2}+\frac{1\mathtt{a}_1}{2\beta_2}\right)Y_4^I Y^J_4\Biggr),
\end{split}
\end{equation}
choosing  $\beta_2=\frac{2+\sqrt{2\sqrt{5}+10}}{1+\sqrt{5}} $ we guarantee that 
\begin{equation}\label{h2B}
\frac{1}{2}+\mathtt{a}_1\frac{1-\sqrt{5}}{2}+ \frac{\mathtt{a}_1\beta_2}{2} =   \frac{1}{2}+\frac{\mathtt{a}_1}{2\beta_2}>0  \quad \text{for all} \ t \in (0, t_0).
\end{equation}
Then we define 
\begin{equation}\label{hB3}
\gamma_1=  \left(\frac{1}{2}+\frac{\mathtt{a}_1}{2\beta_2}\right)\Bigl|_{t=0}=1-\frac{1}{4} \sqrt{2 \sqrt{5}+10},
\end{equation}
we conclude, with the help of   \eqref{h1B}, \eqref{h2B}, \eqref{hB3},  that
\begin{equation}\label{hB4}
h(V,\mathcal{B}_*V) \geq \gamma_1h(V,B^0V). 
 \end{equation} 
Moreover,  using equations,  \eqref{BMEx2},  \eqref{BC7},  \eqref{hB4},  and choosing $m$ such that $\sigma$ is sufficiently small,  it is not difficult to verify that 
\begin{equation}\label{hb5}
 h(V,\tilde{\mathcal{B}}V) \geq \Bigl(\gamma_1-\sigma\Bigr)h(V,B^0 V), 
 \end{equation} 
on $(0, t_0)\times T^1_{\alpha}\times \mathbb{S}^2$, where $t_0$ is given by \eqref{time0}.  Thus we conclude that the existence  of solutions to the conformal wave equations \eqref{eweq} on \eqref{M0} can be obtained from solving the initial value problem 
\begin{align} 
B^{0}\del{t}V^K+\frac{1}{t}\frac{\chi(1-\rho^m)\rho}{m}B^{1} \del{\rho}V^K + \frac{\chi}{t^{\frac{1}{2}}}B^\Sigma \nablasl{\Sigma}V^K &= \frac{1}{t}\tilde{\Bc} \Pbb V^K + \tilde{\Cc} V^K+ \Fc^K  && \text{in} \ \bigl(0, t_0\bigr) \times \mathcal{S}  \label{MwaveJ.1}\\
V^K &= \mathring{V}{}^K && \text{in} \  \bigl\{t_0\bigr\} \times \mathcal{S},  \label{MwaveJ.2}
\end{align} 
for initial data $\mathring{V}{}^K=(\mathring{V}_\Ic^K)$ satisfying the constraints 
\begin{align}
\frac{(1-\rho^m)\rho}{m}\del{\rho}\mathring{V}_4=& -\frac{\sqrt{2}(1+\sqrt{5})}{7+\sqrt{5})}\mathring{V}_0+\frac{3(1+\sqrt{5})}{7+\sqrt{5}}\mathring{V}_1+\rho^m\mathring{V}_4  \quad  \text{in} \quad   \Sigma, \label{ind1}\\
\intertext{and}
\del{\Lambda}\mathring{V}_4=&\frac{\rho^{2m}+2\rho^m-1}{\sqrt{11}\rho^m(1-\rho^m)}\mathring{V}_{\Lambda} \quad  \text{in} \quad   \Sigma.\label{ind2}
\end{align}
The solutions to the equation\eqref{eweq} are determined from the IVP \eqref{MwaveJ.1}, \eqref{MwaveJ.2}, \eqref{ind1}, \eqref{ind2} via
\begin{equation}\label{sol}
u^K (t,r,\theta,\phi)=\frac{(1-r)}{t^{\frac{1}{2}}}V^K_4(t,r^{\frac{1}{m}}, \theta,\phi).
\end{equation}
Using \eqref{sol} we can determine the solution to the system of semi-linear wave equations  presented in \eqref{Mbwave} defined on \eqref{Mregion} using the formula \eqref{uom} which yields
\begin{equation}\label{ubarsl}
\ub^K(\tb,\rb,\thetab, \phib)=\Biggl(\frac{(\tb-\rb)(\tb^2-\rb^2)^{\frac{1}{2}})}{1+\tb-\rb}\Biggr)V^K_4\Bigl(\frac{1}{\tb^2-\rb^2},  \ \frac{1}{(1+\tb-\rb)^{\frac{1}{m}}}, \thetab, \phib\Bigr).
\end{equation}
\subsection{\bf Initial data transformations}
Consider semi-linear equation \eqref{Mbwave} with initial data prescribed on \eqref{indatS0}
\begin{equation}\label{ind3}
(\ub^K, \del{\tb}\ub^K, \del{\rb}\ub^K) = (\vb^K, \wb^K, \zb^K)  \quad \text{in} \quad \bar{\Sigma}_{0},
 \end{equation} 
and the corresponding initial data for the system \eqref{syst1}
\begin{equation}\label{ind4}
(u^K, \del{t}u^K, \del{r}u^K)=(v^K, w^K, z^K)  \quad \text{in} \quad \Sigma.
\end{equation}
Note that  the initial   data \eqref{ind3} is conformally related to \eqref{ind4} as follows
\begin{align}\label{indwd}
v^K(r, \theta, \phi)&=\frac{r^2+2r-1}{2r(1-r)}\vb^K,\\
\label{indw}
w^K(r, \theta, \phi)&=-\frac{2r}{1-r}\Biggl(\frac{r^2+2r-1}{2r(1-r)}\Bigl(\wb^K+\zb \Bigr) +\vb^K  \Biggr), \\
 \label{indz}
z^K(r, \theta, \phi)&=\frac{(r^2+2r-1)^2}{4(1-r)^3r^3}\wb^K+\frac{3r^2-2r+1}{2(1-r)^2r^2}\vb^K+\frac{4r^4-(1-r)^4}{4(1-r)^3r^3}\zb^K.
\end{align}
Using equations \eqref{chv1}, \eqref{chvm}, \eqref{chvm2}, \eqref{rtr} and \eqref{indwd}, \eqref{indw}, \eqref{indz} we obtain the following initial data for the system \eqref{shweq2}
\begin{equation}\label{indv0}
\hat{V}=\begin{pmatrix}
\frac{3}{4(1-\rho^m)}w^K + \rho^mz^K\\
\frac{1}{2\sqrt{2}(1-\rho^m)}w^K+\frac{1+\sqrt{5}}{2\sqrt{2}}\rho^m z^K\\
\frac{\sqrt{11}\rho^m}{\sqrt{2}(r^{2m}+2r^{m}-1)}\del{\theta}v^K\\
\frac{\sqrt{11}\rho^m}{\sqrt{2}(r^{2m}+2r^{m}-1)}\del{\phi}v^K\\
\frac{1}{\sqrt{2}(1-\rho^{m})}v^K
\end{pmatrix}  \quad \quad \text{in} \quad \quad   \Sigma,
\end{equation}
where $v^K=v^K(\rho, \theta, \phi)$,  $w^K=w^K(\rho, \theta, \phi)$,  $z^K=z^K(\rho, \theta, \phi)$.  It is clear that the initial data \eqref{indv0} can be extended for the system \eqref{MwaveJ.1}, \eqref{MwaveJ.2} defined on \eqref{extst2},  and equation  \eqref{indv0} satisfies  the constraints \eqref{ind1}, \eqref{ind2}.  In other words, we can choose initial data $\mathring{V}$ on $\mathcal{S}$ that when restricted to $\Sigma_{\rho_0}$ we obtain 
\begin{equation*}
\mathring{V}\big|_{\Sigma_{\rho_0}} =\hat{V}.
 \end{equation*} 
\section{\bf Construction  of the complete Fuchsian system}\label{complete}

\subsection{The differentiated system}\label{differsys}
Following \cite{Oliynyk2021107} we differentiate the system \eqref{MwaveE} by applying the Levi-Civita connection $\mathcal{D}_j$ of the Riemannian metric on $\mathcal{S}$ 
\begin{equation*}
q=q_{ij}dx^i \otimes dx^j \coloneqq  d\rho \times d\rho + \gsl,
\end{equation*}  
note that in the coordinates $x^{i}=(\rho, \theta , \phi)$,  the operator  $\mathcal{D}$ takes the form  
\begin{equation}\label{derD}
\mathcal{D}_j =\delta^1_j \del{\rho} +\delta^{\Lambda}_j \nablasl{\Lambda},
\end{equation}
where $\nablasl{\Lambda}$ is the Levi-Civita connection associated to the metric $\gsl_{\Lambda\Sigma}$ on $\mathbb{S}^2$. Applying $\mathcal{D}_j$ to \eqref{MwaveE} and multiplying by $t^{\kappa}$, where $\kappa$ is a small positive constant  that we fix below,  we obtain  
\begin{equation}\label{difs}
B^0 \del{t}W_j^K+\frac{\chi (1-\rho^m)\rho}{tm}B^1\del{\rho}W_j^K + \frac{\chi}{t^{\frac{1}{2}}}B^{\Sigma}\nablasl{\Sigma}W_j^K = \frac{1}{t}\Bigl(\tilde{\mathcal{B}}\mathbb{P}+\kappa B^0 \Bigr)W_j^K+ \frac{1}{t}\mathcal{Q}_j^K +\mathcal{H}^K_j,
\end{equation}
where
\begin{equation}\label{Wvar}
W^K_j=(W^K_{j\mathcal{I}})\coloneqq (t^{\kappa}D_jV^K_{\mathcal{I}}),
\end{equation}
\begin{equation} \label{Qcdef}
 \Qc_j^K = - \frac{\rho^{3m}}{2(\rho^{2m}-(1-\rho^m)^2t)}t^{\kappa}b^{K}_{IJ}\chi(\rho) \mathcal{D}_j(V^I_0 V_0^J),
\end{equation}
and
\begin{align} 
\Hc^K_j = & - \frac{1}{t}\mathcal{D}_j\biggl(\frac{\chi(1-\rho^m)\rho}{m}B^1\biggr)W^K_1-\frac{1}{t^{\frac{1}{2}}}\mathcal{D}_j\Bigl(\chi B^{\Sigma} \Bigr)W^K_{\Sigma}+ t^{\kappa-\frac{1}{2}} B^\Sigma[\nablasl{\Sigma},\Dc_j]V^K +t^{\kappa-1}\mathcal{D}_j(\tilde{\Bc})\mathbb{P}V^K  \notag \\
&+t^{\kappa-\frac{1}{2}}\mathcal{D}_j(\tilde{\Cc})V^K+t^{-\frac{1}{2}}\tilde{\Cc} W^K_j+ t^{\kappa-1}\Dc_j\Bigl(\frac{\rho^{3m}}{\rho^{2m}-(1-\rho^m)^2t}\bb^{K}_{IJ}\chi(\rho)\Bigr) V^I_0 V_0^J \ev_0+t^{\kappa}\Dc_j \Gc.\label{HcKdef}
\end{align}

\subsection{The asymptotic equation}
Using the notation introduced in \eqref{Wvar},  consider now the the first equation of the extended system \eqref{MwaveE}
\begin{equation} \label{MwaveL}
\begin{split}
&\del{t}V_0^K +\frac{\chi(1-\rho^m)\rho}{t^{\kappa}m}\Bigl(\mathtt{a}_0W^K_{10}+t^{\frac{1}{2}}\mathtt{a}_1W^K_{11}\Bigr) -\frac{1}{t^{\frac{1}{2}+\kappa}}\frac{(1-\rho^m)(1+t)}{\mathtt{p}(\rho^{2m}-(1-\rho^m)^2t)}\gsl^{\Lambda\Sigma}W^K_{\Lambda\Sigma}= \mathtt{a}_0(2+\chi(3+\rho^m))V_0^K+\\
&\frac{\mathtt{a}_1}{t^{\frac{1}{2}}}(2+\chi(3+\rho^m))V^K_1- \frac{\rho^{3m}}{2t(\rho^{2m}-(1-\rho^m)^2t)}b^{K}_{IJ}\chi(\rho) V^I_0 V_0^J + \Gc_{0}^K,
\end{split}
\end{equation}
From the definition of the asymptotic equation introduced in \cite{Oliynyk2021107} and the equations \eqref{asympprop1.1},\eqref{asympprop1.2},  \eqref{MwaveL},   we see that the asymptotic equation associated to the system \eqref{MwaveE}, is given by
\begin{equation}\label{aeqi2}
\del{t}\xi= \frac{1}{t}Q(\xi),
\end{equation}
 where
\begin{equation}\label{aeqi3}
Q(\xi)=Q(\xi^K)=- \frac{\rho^{3m}}{2(\rho^{2m}-(1-\rho^m)^2t)}b^{K}_{IJ}\chi(\rho) \xi^I \xi^J.
\end{equation}
The asymptotic equation \eqref{aeqi2} involves the most singular term in the quadratic  nonlinearities.  This term is challenging due to its singularity at $t=0$.    Thus the  system \eqref{MwaveE} is not yet  in the form required in \cite{BOOS:2020}.   In order to remove the singular term $\frac{1}{t}Q^K$,  we define  the flow\footnote{Note that the flow depends on $y=(y^i)=(\rho,\theta,\phi)\in \Sc$ through the coefficients $\chi \rho^m \bb^K_{IJ}$, which are smooth functions on $\Sc$.} $\Fsc(t,t_0,y, \xi)=(\Fsc^K(t,t_0,y,\xi))$ such that
\begin{align}
\del{t}\Fsc(t,t_0,y,\xi) &=  \frac{1}{t}Q\bigl(\Fsc(t,t_0,y,\xi)\bigr), \label{asympIVP.1} \\
\Fsc(t,t_0,y,\xi)&= \mathring{\xi}.\label{asympIVP.2}
\end{align}
where $t_0$ is defined in \eqref{time0}.  For fixed $(t,t_0,y)$, the flow $\Fsc(t,t_0,y, \xi)$ maps $\Rbb^{N}$ to itself, and consequently, the derivative $D_\xi F(t,t_0,y,\xi)$ 
defines a linear map from $\Rbb^{N}$ to itself, or equivalently, a $N \times N$-matrix.  Using the flow  we define a new set of variables $Y(t,y)=(Y^K(t,y))$ via
\begin{equation}\label{Ydef}
V_0(t,y) = \Fsc(t,t_0,y,Y(t,y))
\end{equation}
where
\begin{equation}\label{V0def}
V_0= (V^K_0).
\end{equation}
Using \eqref{Ydef} and  equations  \eqref{MwaveL},  \eqref{asympIVP.1} it is straightforward to verify  that  $Y$ satisfies
\begin{equation}  \label{MwaveM}
\del{t}Y= \Lsc\Gsc
\end{equation}
where
\begin{equation} \label{Lscdef}
 \Lsc =(D_{Y}\Fsc(t, t_0, y, Y))^{-1}
 \end{equation}
 and
 \begin{equation}\label{LscdeG}
 \begin{split}
\Gsc = &-\frac{\chi(1-\rho^m)\rho}{t^{\kappa}m}\Bigl(\mathtt{a}_0W^K_{10}+t^{\frac{1}{2}}\mathtt{a}_1W^K_{11}\Bigr) +\frac{1}{t^{\frac{1}{2}+\kappa}}\frac{(1-\rho^m)(1+t)}{\mathtt{p}(\rho^{2m}-(1-\rho^m)^2t)}\gsl^{\Lambda\Sigma}W^K_{\Lambda\Sigma}+ \mathtt{a}_0(2+\chi(3+\rho^m))V_0^K+\\
&\frac{\mathtt{a}_1}{t^{\frac{1}{2}}}(2+\chi(3+\rho^m))V^K_1 + \Gc_{0}^K.
\end{split}
\end{equation}

\subsubsection{\bf Asymptotic flow assumptions\label{Asymptoticflowassump}}
The flow \eqref{asympIVP.1}, \eqref{asympIVP.1} in conjunction with equations \eqref{Ydef}, \eqref{V0def}, \eqref{MwaveM}, \eqref{Lscdef}, \eqref{LscdeG} is of the same form as the flow in section  3.4 of \cite{Oliynyk2021107}.  Therefore we say that the flow $\Fsc(t,t_0,y,\xi)=(\Fsc^K(t,t_0,y,\xi))$ satisfies the following: Given any $\Ntt\in \Zbb_{\geq 0}$,  there exist constants $R_0>0$, $\ep \in (0,1)$ and $C_{k\ell }>0$,  where $k,\ell\in \Zbb_{\geq 0}$ and $0 \leq k+\ell\leq \Ntt$, and a function $\omega(R)$ satisfying $\lim_{R\searrow 0} \omega(R)=0$
such that
\begin{gather} 
\bigl|\Fsc(t,t_0,y,\xi)\big|\leq \omega(R),  \label{flowassump.1}
\intertext{and}
\Bigl| D_\xi^k \Dc^\ell \Fsc\bigl(t,t_0,y,\xi\bigr)\Big| 
+\Big|D_\xi^k \Dc^\ell  \bigl(D_\xi\Fsc\bigl(t,t_0,y,\xi\bigr)\bigr)^{-1}\Big| 
\leq \frac{1}{t^{\ep}}C_{k\ell},  \label{flowassump.2}
\end{gather}
for all $(t,y,\xi) \in (0,t_0]\times \Sc \times B_R(\Rbb^{N})$ and $R\in (0,R_0]$.  From here,  one sees that the maps $\Ftt$ and $\check{\Ftt}$ 
defined by
\begin{equation} \label{flowassump1}
\Ftt(t,y,\xi)=t^\ep \Fsc\bigl(t, t_0,y,\xi\bigr) \AND = \check{\Ftt}(t,y,\xi)=t^\ep \Bigl(D_\xi\Fsc\bigl(t, t_0,y,\xi\bigr)\Bigr)^{-1},
\end{equation}
satisfy
$\Ftt \in C^0\Bigl(\bigl[0, t_0\bigr],C^{\Ntt}(\Sc\times B_R(\Rbb^N),\Rbb)\Bigr)$ and 
$\check{\Ftt} \in C^0\Bigl(\bigl[0, t_0\bigr],C^{\Ntt}(\Sc\times B_R(\Rbb^N),\Mbb{N})\Bigr)$.
Moreover,  since $\xi=0$  solves the asymptotic equation \eqref{aeqi2},  the flow  satisfies $\Fsc(t,t_0,y,0)=0$,  thus
\begin{equation} \label{flowassump2}
\Ftt(t,y,0)=0 
\end{equation}
for all $(t,y)\in [0, t_0]\times \Sc$.
\begin{prop} \label{asympprop}
Suppose the bounded weak null condition holds (see Definition \ref{bwnc}). Then there exists a $R_0\in (0,\Rc_0)$
such that the flow $\Fsc(t,t_0,y,\xi)$ of the asymptotic equation \eqref{aeqi2} satisfies the flow assumptions
\eqref{flowassump.1}-\eqref{flowassump.2} for this choice of $R_0$ and any choice of $\ep \in (0,1)$.
\end{prop}
\begin{lem} \label{asymplem}
For any $R\in (0,\Rc_0]$, the solutions $\xi$ of the asymptotic IVP \eqref{asympprop1.1}-\eqref{asympprop1.2} exist
for $t\in (0, t_0]$ and satisfies
\begin{equation}
\sup_{0<t\leq t_0}|\xi(t)| \leq \frac{C}{\Rc_0} R \label{asymplem0}
\end{equation}
for any choice of initial data that is bounded by $|\mathring{\xi}| \leq R$.
\end{lem}
The proof for Proposition \eqref{asympprop} and Lemma \eqref{asymplem} follow directly  from Proposition 3.2 from \cite{Oliynyk2021107}.  Note that  the solution $\xi=(\xi^K)$ will  depend implicitly on $y\in \Sc$ and the initial data $\mathring{\xi}$.  Note that for a fixed  $\ep >0$ one can differentiate  the asymptotic equation \eqref{aeqi2} with respect to $y=(y^i)$ and define  
\begin{equation}\label{asympprop3}
\eta^K_i= t^{\ep}\Dc_i\xi^K,
\end{equation}
which satisfies the  equation
\begin{equation} \label{asympprop3a}
\del{t} \eta^K_i =  \frac{1}{t}\bigl( \ep \delta^K_J  - \frac{\rho^{3m}\chi}{2(\rho^{2m}-(1-\rho^m)^2t)}\bigl(\bb^K_{JI}
+ \bb^K_{IJ}\bigr)\xi^I \bigr) \eta^J_i - \frac{1}{t^{1-\ep}}\Dc_i\Bigl( \frac{\rho^{3m}\chi \bb^K_{IJ}}{2(\rho^{2m}-(1-\rho^m)^2t)}\Bigr) \xi^I \xi^J.   
\end{equation}
By contracting \eqref{asympprop3a} with $\delta_{LK}\delta^{ki}\eta_k^L$ and  defining the Euclidean norm $|\eta| = \sqrt{\delta_{KL}\delta^{ij}\eta^K_i \eta^L_j}$,   where $\eta=(\eta_i^K)$,  we obtain 
\begin{equation} \label{asympprop4}
\begin{split}
\frac{1}{2}\del{t}|\eta|^2 = \frac{1}{t}\bigl(\ep |\eta|^2 - \frac{\rho^{3m}\chi}{2(\rho^{2m}-(1-\rho^m)^2t)}\bigl(\bb^K_{JI}
+ \bb^K_{IJ}\bigr) \delta_{LK}\xi^I  \delta^{ki}\eta_k^L\eta^J_i \bigr) -\\
\frac{1}{t^{1-\ep}}\delta_{LK}\delta^{ki}\eta_k^L\Dc_i\Bigl(- \frac{\rho^{3m \chi \bb^K_{IJ}}\chi}{2(\rho^{2m}-(1-\rho^m)^2t)}\Bigr) \xi^I \xi^J.
\end{split}
\end{equation}
Notice that $\frac{\rho^{3m}\chi}{2(\rho^{2m}-(1-\rho^m)^2t)}\bb^K_{JI}$ is a smooth  bounded function  on $\Sc$,  and so are their derivatives.  Using this fact and the bound on $\xi$ given by  \eqref{asymplem0},  we deduce from  \eqref{asympprop4}
and the Cauchy Schwartz inequality that there exist $\upsilon>0$ such that for any $\upsilon \in (0,\ep)$ there exists constants $R_0\in (0,\Rc_0]$ and $C>0$ such that the energy inequality
\begin{equation*}
\frac{1}{2}\del{t}|\eta|^2 \geq \frac{1}{t}(\ep-\upsilon) |\eta|^2 -\frac{C}{t^{1-\ep}}|\eta|
\end{equation*}
holds for any given $R\in (0,R_0]$ and for all $t\in (0,t_0]$.  Which leads to 
\begin{equation*}
\del{t}|\eta| \geq \frac{1}{t}(\ep-\upsilon) |\eta| -\frac{C}{t^{1-\ep}}.
\end{equation*}
Applying of Gr\"{o}nwall's inequality\footnote{Here, we are using the following form of Gr\" {o}nwall's inequality: if $x(t)$ satisfies
$x'(t) \geq a(t)x(t) -h(t)$, $0<t\leq t_0$, then $x(t)\leq x(t_0)e^{-A(t)}+ \int^{t_0}_t e^{-A(t)+A(\tau)}h(\tau) \, d\tau$ where
$A(t)= \int^{t_0}_t a(\tau)\, d\tau$. In particular, we observe from this that if, $x(t_0)\geq 0$ and  $a(t)=\frac{\lambda}{t}-b(t)$, where  $\lambda \in \Rbb$ and $\bigl|\int^{t_0}_t b(\tau) \, d\tau \bigr|\leq r$, then
\begin{equation*}
x(t) \leq e^{r}x(t_0)\left(\frac{t}{t_0}\right)^\lambda +   e^{2r}t^\lambda \int^{t_0}_t \frac{|h(\tau)|}{\tau^\lambda}\, d\tau
\end{equation*}
for $0\leq t < t_0$.
} we obtain 
\begin{equation} \label{asympprop4a}
|\eta(t)| \leq |\eta(t_0)| t^{\ep-\upsilon}+   t^{\ep-\upsilon}\int^{t_0}_t \frac{C}{t^{1-\upsilon}}\, d\tau = 
t^{\ep-\upsilon}|\eta(t_0)| +\frac{1}{\upsilon}t^{\ep-\upsilon}\bigl(t_0-t^\upsilon\bigr).
\end{equation}
Using  inequality \eqref{asympprop4a},  \eqref{asymplem0},  definition \eqref{asympprop3} and since $\xi(t)=\Fsc(t, t_0,y,\mathring{\xi})$, we conclude that there exist constants $C_0, C_{01}>0$ such that the flow $\Fsc$ satisfies the bounds
\begin{equation*}  
|\Fsc(t, t_0, y,  \mathring{\xi})| \leq C_0 R 
\AND
|\Dc\Fsc(t, t_0, y,\mathring{\xi})| \leq \frac{1}{t^{\upsilon}}C_{01} 
\end{equation*}
for all $(t,y,\mathring{\xi})\in (0, t_0]\times \Sc \times B_{R}(\Rbb^N)$, $R\in (0,R_0]$.
Using similar arguments it is not difficult to show that 
\begin{equation} \label{asympprop5}
\del{t} D_{\mathring{\xi}}\xi =  \frac{1}{t}L D_{\mathring{\xi}}\xi
\end{equation}
where 
\begin{equation*}
L=(L^K_J):=\frac{\rho^{3m}\chi}{2(\rho^{2m}-(1-\rho^m)^2t)}\bigl(\bb^K_{JI}
+ \bb^K_{IJ}\bigr).
\end{equation*}
Multiplying \eqref{asympprop5} on the right by  $(D_{\mathring{\xi}}\xi)^{-1}$ leads to
\begin{equation} \label{asympprop6}
\del{t} ( (D_{\mathring{\xi}}\xi)^{-1})^{\tr} =  -\frac{1}{t} 
L^{\tr} ((D_{\mathring{\xi}}\xi)^{-1})^{\tr}
\end{equation}
Now using  \eqref{asympprop5} and \eqref{asympprop6} and multiplying by 
$t^\ep$,  we get
\begin{align*}
\del{t}( t^\ep D_{\mathring{\xi}}\xi) &=  \frac{1}{t}
\bigl((\ep + L\bigr) t^\ep D_{\mathring{\xi}}\xi  
\intertext{and}
\del{t} (t^\ep (D_{\mathring{\xi}}\xi)^{-1})^{\tr} &=  \frac{1}{t}
\bigl(\ep-L^{\tr}\bigr)(t^\ep (D_{\mathring{\xi}}\xi)^{-1})^{\tr}.
\end{align*}
These equations are similar to  \eqref{asympprop3a}, and thus we can use the same reasoning  to derive  estimates for $t^\ep D_{\mathring{\xi}}\xi$ and $(t^\ep(D_{\mathring{\xi}}\xi)^{-1})^{\tr}$.  Therefore  we conclude  that there exist a constant $C_{10}>0$ such that  
\begin{equation*}
\bigl|D_{\mathring{\xi}}\xi\bigr|+ \bigl|(D_{\mathring{\xi}}\xi)^{-1}\bigr|\leq \frac{1}{t^{\upsilon}}C_{10}
\end{equation*}
holds for $0<t\leq t_0$.  This estimate leads to 
\begin{equation*} 
\bigl|D_{\mathring{\xi}}\Fsc(t, t_0,y,\mathring{\xi})\bigr|
+\bigl|\bigl(D_{\mathring{\xi}}\Fsc(t, t_0,y,\mathring{\xi})\bigr)^{-1}\bigr|\leq \frac{1}{t^{\upsilon}}C_{10},
\end{equation*}
for all $(t,y,\mathring{\xi})\in (0, t_0]\times \Sc \times B_{R}(\Rbb^N)$ and $R\in (0,R_0]$.  We can also use similar arguments  to derive, for any fixed $\Ntt \in \Zbb_{\geq 1}$,  the bounds
\begin{equation*}
\bigl|D^k_{\mathring{\xi}}\Dc^\ell\xi\bigr| +\bigl|D^k_{\mathring{\xi}}\Dc^\ell \bigl(D_{\mathring{\xi}}\xi)^{-1}\bigr|\leq \frac{1}{t^\upsilon}C_{kl} 
\end{equation*}
on the higher derivatives for $1\leq k + \ell \leq \Ntt$.  Therefore  the flow satisfies the bounds
\begin{equation*}
|D_{\mathring{\xi}}^k \Dc^\ell \Fsc(t, t_0,y,\mathring{\xi})| \leq \frac{1}{t^{\upsilon}}C_{\ell k},
\end{equation*}
hold for all $(t,y,\mathring{\xi})\in (0, t_0]\times \Sc \times B_{R}(\Rbb^N)$, $2\leq k + \ell \leq \Ntt$, and $R\in (0,R_0]$. 
\subsubsection{\bf The complementary  variable $X^K$ }
We complement equations \eqref{difs} and \eqref{MwaveM} with a third system obtained from a  rescaling of  the projection operator $\mathbb{P}$ applied to the variable \eqref{MwaveE},  we defined the variable $X^K $ by
\begin{equation}\label{XKdef}
X^K = \frac{1}{t^\nu}\Pbb V^K,
\end{equation}
where $\nu\geq 0$ is a constant that we fix below.  Using equations  \eqref{Pop},  \eqref{Poppr1}, \eqref{FcKdef},  \eqref{Gcdef},  \eqref{Wvar} it is not difficult to see that  $X^K$ satisfies
\begin{align}
B^0\del{t}X^K + \frac{1}{t} \frac{\chi (1-\rho^m)\rho}{m}B^1\del{\rho}X^K  +\frac{\chi}{t^{\frac{1}{2}+\nu \kappa}}B^{\Lambda}\nablasl{\Lambda}X^K&= \frac{1}{t}(\tilde{\Bc} - \nu B^0) X^K +\Kc^K
\label{MwaveN}
\end{align}
where
\begin{equation}\label{KcKdef}
\Kc^K=-\frac{\chi(1-\rho^m)\rho}{mt^{1+\kappa+\nu}}[\mathbb{P}, B^1]W^K_1- \frac{\chi}{t^{\frac{1}{2}+\kappa+\nu}}[\Pbb,  B^{\Lambda}] W_{\Lambda}^K +\Pbb\Cc \biggl(\frac{1}{t^{\nu}}\Pbb^\perp V^K + X^K\biggr) + \frac{1}{t^\nu}\Pbb \Gc_0^K + \frac{1}{t^{\frac{1}{2}+\nu}}\Pbb \Gc^K_1 
\end{equation}
and the projection operator and complement operator $\Pbb,  \Pbb^\perp,$ satisfy \eqref{Poppr1}.
\subsubsection{\bf The complete Fuchsian system}
Now  using the variables \eqref{difs}, \eqref{MwaveM} and \eqref{MwaveN} we can write the complete Fuchsian system  as follows
\begin{align}
 A^0\del{t}Z + \frac{1}{t}\frac{\chi (1-\rho^m)\rho}{m}A^1 \del{\rho}Z + \frac{\chi}{t^{\frac{1}{2}}}
A^\Sigma \nablasl{\Sigma}Z
&= \frac{1}{t}\Ac \Pi Z +  \frac{1}{t}\Qc + \Jc \label{MwaveO}
\end{align}
where
\begin{align}
Z &=  \begin{pmatrix} W^K_j & X^K & Y \end{pmatrix}^{\tr} , \notag\\
A^0 &= \begin{pmatrix} B^0 & 0 & 0 \\ 0 & B^0 & 0 \\ 0 & 0 & \id \end{pmatrix}, \label{A0def} \\
A^1 & = \begin{pmatrix}  B^1\delta^j_k  & 0 & 0 \\ 0 & B^1 & 0 \\ 0 & 0 &0 \end{pmatrix}, \label{A1def} \\
A^\Sigma & = \begin{pmatrix}  B^\Sigma  & 0 & 0 \\ 0 & B^\Sigma  & 0 \\ 0 & 0 & 0  \end{pmatrix}, \label{ASigmadef} \\
\Ac & = \begin{pmatrix} \tilde{\Bc}\Pbb + \kappa B^0 & 0 & 0 \\ 0 & \tilde{\Bc}-\nu B^0 & 0 \\ 0 & 0 & \id  \end{pmatrix}, \label{Acdef} \\
\Pi & = \begin{pmatrix}  \id & 0 & 0 \\ 0 & \id & 0 \\ 0 & 0 & 0  \end{pmatrix}, \label{Pidef}\\
\Qc &= \begin{pmatrix} \Qc^K_j & 0 & 0 \end{pmatrix}^{\tr} \notag
\intertext{and}
\Jc &= \begin{pmatrix}\Hc^K_j &  \Kc^K & \Lsc\Gsc \end{pmatrix}^{\tr}.\label{Jcdef}
\end{align}
where $Z$ is a time-dependent section of the vector bundle
\begin{equation*}
\Wbb = \bigcup_{y\in \Sc} \Wbb_{y}
\end{equation*}
over $\Sc$ with fibers
$\Wbb_y =\Bigl(T^*_{y}\Sc \times T^*_{y}\Sc \times \bigl(\text{T}^*_{y}\Sc \otimes T^*_{\text{pr}(y)}\mathbb{S}^2\bigr)\times T^*_{y}\Sc\Bigr)^N\times \Vbb_y^N\times \Rbb^N$ and  $\text{pr} :  \Sc  \longrightarrow \mathbb{S}^2$ is the canonical projection and $\Vbb_y=\Rbb\times \Rbb\times \text{T}^*_{\text{pr}(y)}\mathbb{S}^2\times\Rbb$.  Taking  $\grave{Z},  Z \in \mathbb{W}$,  with $ \grave{Z}= ( \grave{W}{}^K_j, \grave{X}{}^K , \grave{Y}),$ we introduce the  inner-product on $\Wbb$ via
\begin{equation} \label{hcdef}
\hc(Z,\grave{Z}) = \delta_{KL}q^{ij}h(W^K_i, \grave{W}{}^L_j)
+ \delta_{KL}h(X^K,\grave{X}{}^L)+\delta_{KL}Y^K \grave{Y}^L,
\end{equation} 
where $h(\cdot,\cdot)$ is the inner-product defined previously by \eqref{hdef}.  It is clear from this definition of inner product and equations \eqref{A0def}-\eqref{ASigmadef} that $A^0$, $A^1, A^\Sigma \eta_\Sigma$ and $\Pi$,  are symmetric with respect to the inner-product \eqref{hcdef}. Using  the connection $\Dc_j$ defined in \eqref{derD} we can  verify that the inner-product \eqref{hcdef} is compatible with the connection $\Dc_j$ defined in \eqref{derD},  that is   $\Dc_j\bigl(\hc(Z,\grave{Z})\bigr)= \hc(\Dc_j Z,\grave{Z})+\hc(Z, \Dc_j \grave{Z})$.   The operator $\Pi$ defined in  \eqref{Pidef} is a  projection operator,  together with its  complementary projection operator  $\Pi^\perp$  they satisfy 
\begin{equation*}
\Pi^2 = \Pi , \quad \Pi^\perp = \id -\Pi.
\end{equation*}
Moreover,   from the definitions \eqref{A0def}, \eqref{A1def} and \eqref{Acdef}  we see that
\begin{equation*}
[A^0,\Pi]=[\Ac,\Pi]=0, \quad \quad [\Pi^\perp, A^i] =0,  \quad     A^i \Pi^\perp =0.
\end{equation*}
and
\begin{equation*} 
[\Pi, A^i] =0 ,  \quad A^i \Pi =  A^i ,
\end{equation*}
Next, we see from \eqref{InpV}, \eqref{A0def} and \eqref{hcdef} that $A^0$ satisfies
\begin{align*}
\hc(Z,A^0 Z) &= \delta_{KL}q^{ij}h(W^K_i,B^0 W^L_{j})
+ \delta_{KL}h(X^K,B^0 X^L)+\delta_{KL}Y^K Y^L \\
 &= \delta_{KL}q^{ij}h(W^K_i,W^L_j)
+ \delta_{KL}h(X^K,X^L)+\delta_{KL}Y^K Y^L,
\end{align*}
and hence, that
\begin{equation*} 
 \hc(Z,A^0 Z)=\hc(Z,Z).
\end{equation*}
Then using equations  \eqref{Poppr1},  \eqref{Poppr2}, \eqref{hdef}, \eqref{hb5}, \eqref{A0def}, \eqref{Acdef} and \eqref{hcdef}
we can show that 
\begin{equation*}
\hc(Z,\Ac Z) \geq \kappa \hc(Z,A^0 Z)
\end{equation*} 
as long as the constants   $\kappa>0,  \nu>0$ satisfy   $\gamma_1-\sigma \geq \kappa+\nu $.  
Using the inequality \eqref{boundpB1} and definition \eqref{A1def}, it is not difficult to show that there exist an integer $m \geq 1$ and constant $\sigma_2 > 0$ such that  
\begin{equation} \label{A1bnd}
\biggl|\del{\rho}\biggl(\frac{\chi(1-\rho^m)\rho}{m}B^1\biggr)\biggr| + \biggl|\del{\rho}\biggl(\frac{\chi(1-\rho^m)\rho}{m} A^1\biggr)\biggr| < \sigma_2 \quad \text{in $(0,t_0)\times \Sc$}.
\end{equation}

\subsubsection{ \bf The source term $\Jc$:}
It is a straightforward calculation to verify that the source term defined in \eqref{Jcdef} is of the same form as the source term in section 3.6.2 from \cite{Oliynyk2021107},  and satisfies the same conditions.  Using \eqref{VKvar},  \eqref{Pop},  \eqref{Poppr1} notice that we can decompose the variable $V^K$ as follows 
\begin{equation} \label{VKdecomp1}
V^K(t,y) = \Pbb V^K(t,y) + \Pbb^\perp V^K(t,y),
\end{equation}
and using definition \eqref{XKdef} we see that 
\begin{equation}
\Pbb V^K(t,y) = t^\nu X^K(t,y).  \label{VKdecomp2}
\end{equation}
similarly,  using \eqref{VKvar},  \eqref{Pop},  \eqref{Poppr1},  \eqref{Ydef}, \eqref{flowassump.1} we see that
\begin{equation}
 \Pbb^\perp V^K(t,y) =  \frac{1}{t^{\ep}}( t^\ep V_0^K(t,y))\ev_0 =   \frac{1}{t^{\ep}}\Ftt^K\bigl(t,y,Y(t,y)\bigr)  \bigr)\ev_0,\label{VKdecomp3}
\end{equation}
and from  definition  \eqref{Wvar} we see that the derivative $\Dc_j V^K$  can be written as
\begin{equation} \label{VKdecomp4}
\Dc_j V^K(t,y) = \frac{1}{t^{\kappa}}W_j^K(t,y).
\end{equation}  
Finally, by   using equations  \eqref{Lscdef} and \eqref{flowassump1} we can write  the map $\Lsc$ as follows 
\begin{equation} \label{Lscmap}
\Lsc = \frac{1}{t^\ep}\check{\Ftt}\bigl(t,y,Y(t,y)\bigr).
\end{equation}
Now we can write the component \eqref{KcKdef} of the  source term \eqref{Jcdef},   using  equations  \eqref{VKdecomp1},  \eqref{VKdecomp2},  \eqref{VKdecomp3}, \eqref{VKdecomp4},  \eqref{Lscmap}  as well as  \eqref{GcK0def}-\eqref{GcK1def},  which leads to
\begin{align*}
\Kc^K = - \frac{1}{t^{\frac{1}{2}+\kappa+\nu}}&\Pbb B^\Sigma(t,y) W_\Sigma^K(t,y) +\frac{1}{t^{\nu+\ep}}\Ftt^K\bigl(t,y,Y(t,y) \bigr)
\Pbb\Cc(t)\ev_0 + \Pbb\Cc(t) X^K(t,y) \\
+& \frac{1}{t^{\nu+2\ep}}
 \Pbb \Gc^K_0\Bigl(t^{\frac{1}{2}},t,\chi(\rho)\rho^m,\Ftt\bigl(t,y,Y(t,y)\bigr)\ev_0 ,\Ftt\bigl(t,y,Y(t,y)\bigr)\ev_0 \Bigr) \\
 +&  \sum_{a=0}^1\biggl\{ \frac{1}{t^{\frac{a}{2}+\ep}}\biggr[
 \Pbb \Gc^K_a\Bigl(t^{\frac{1}{2}},t,\chi(\rho)\rho^m,\Ftt\bigl(t,y,Y(t,y)\bigr)\ev_0 ,X(t,y)\Bigr) \\
 +& \Pbb \Gc^K_a\Bigl(t^{\frac{1}{2}},t,\chi(\rho)\rho^m,X(t,y),\Ftt\bigl(t,y,Y(t,y)\bigr)\ev_0 \Bigr)   \biggr] +  \frac{1}{t^{\frac{a}{2}-\nu}}\Pbb \Gc^K_a\Bigl(t^{\frac{1}{2}},t,\chi(\rho)\rho^m,X(t,y),X(t,y)\Bigr)   \biggr\}.
\end{align*}
We can write a similar expansion for the components of source term $\Hc^K_j$ and $\Lsc\Gsc$,   in terms of $W^K_j$, $X^K$ and $Y^K$,  by using the definitions\eqref{HcKdef} and \eqref{Lscdef},  \eqref{flowassump2},  \eqref{A1bnd} and using the constants  $\ep,\kappa,\nu\geq 0$.  Thus we conclude that the source term \eqref{Jcdef} can be expanded as  
\begin{equation*}
\begin{split}
\Jc = \biggl(\frac{1}{t^{3\ep}}+\frac{1}{t^{\nu+2\ep}}+\frac{1}{t^{1-\kappa+2\ep}}\biggr)\Jc_0\bigl(&t,y,Z(t,y)\bigr) + \biggl(\frac{1}{t^{\frac{1}{2}+\kappa+\ep}}+\frac{1}{t^{\frac{1}{2}+2\ep-\nu}}\biggr) \Jc_1\bigl(t,y,Z(t,y)\bigr) \notag \\
&+ \frac{1}{t}\bigl(\sigma_2 + t^{\frac{1}{2}-\kappa-\nu}+
t^{\frac{1}{2}-\ep}
+ t^{\frac{1}{2}-\kappa-\ep}+t^{2\nu-\ep}\bigr)\Jc_2\bigl(t,y,Z(t,y)\bigr)
\end{split}
\end{equation*}
where  $\Jc_a\in C^0([0,1], C^\Ntt(\Sc\times B_R(\Wbb),\Wbb))$, $a=0,1,2$, for any fixed $\Ntt\in \Zbb_{\geq 0}$, and these maps satisfy\footnote{Here, we are using are the order notation $\Ord(\cdot)$ from \cite[\S 2.4]{BOOS:2020} where the maps are
finitely rather than infinitely differentiable.}
\begin{gather*}
\Jc_0 = \Ord(Z), \quad \Jc_1 = \Ord(\Pi Z), \quad
\Pi\Jc_2  = \Ord(\Pi Z) \AND
\Pi^\perp\Jc_2  = \Ord(\Pi Z\otimes \Pi Z). 
\end{gather*}

Moreover,  we  can choose the same restrictions for the constants $\kappa,\nu \in \Rbb_{>0}$ to satisfy the inequalities
\begin{equation} \label{kappa-nu-fix1}
\kappa+\nu \leq \gamma_1+\sigma < \frac{1}{2}-\ep, \quad 2\ep < \kappa < 1-\ep,  \quad  \ep < 2\nu 
\end{equation}
These inequalities lead to
\begin{gather*}
3\ep \leq 1-\kappa + 2\ep, \quad \nu+2\ep \leq 1-\kappa + 2\ep, \quad  0<2\nu-\ep, \quad 
0<\frac{1}{2}-\kappa-\ep, \quad  0<\frac{1}{2}-\kappa-\nu,  \\
\frac{1}{2}+2\ep -\nu \leq 1-\frac{\kappa}{2}+\ep, \quad \frac{1}{2}+\kappa+\ep \leq 1-\frac{\kappa}{2}+\ep
\AND 0<\kappa - 2\ep \leq 1,
\end{gather*}
therefore  the map $\Jc_a$,  is of the same form as the source term from section 3.6.2,  which reads as 
 \begin{align*}
\Jc = \frac{1}{t^{1-\kappa+2\ep}}\Jc_0\bigl(t,y,Z(t,y)\bigr)+\frac{1}{t^{1-\frac{\kappa}{2}+\ep}}\Jc_1\bigl(t,y,Z(t,y)\bigr) + \frac{1}{t}\bigl(\sigma_2 + t^{\tilde{\ep}}\bigr)\Jc_2\bigl(t,y,Z(t,y)\bigr)
\end{align*}
for some suitably small constant $\tilde{\ep}>0$.  Note that we can choose our  constant $\epsilon,  \sigma_2 >0$ as small as we like.

\begin{thm} 
Suppose $k\in \Zbb_{\geq 5}$, $\rho_0>0$, there exist $t_0$ such that the extended system \eqref{MwaveJ.1},  \eqref{MwaveJ.2} is symmetric hyperbolic for all $t\in(0,t_0]$,  the asymptotic flow assumptions \eqref{flowassump.1}-\eqref{flowassump.2} are satisfied
for constants $\Ntt \in \Zbb_{\geq k}$,  and the constants $\kappa,\nu, \epsilon,  \in \Rbb_{>0}$ satisfy the inequalities \eqref{kappa-nu-fix1}, and $\zc \in (0,\kappa)$,  then 
\begin{enumerate}
\item  There exist $Z$ such that
\begin{equation*}
Z\in C^0\bigl((0, t_0],H^{k}(\Sc,\Wbb)\bigr)\cap C^1\bigl((0, t_0],H^{k-1}(\Sc,\Wbb)\bigr),
\end{equation*}
which satisfies an energy estimate of the form
\begin{equation*} 
\norm{Z(t)}_{H^k(\Sc)}^2 + \int^{t_0}_{t} \frac{1}{\tau} \norm{\Pi Z(\tau)}^2_{H^k(\Sc)}\, d\tau
\leq C_E^2\norm{Z(t_0)}^2,
\end{equation*}
moreover,
\begin{equation*}
\norm{V_0(t)}_{H^k(\Sc)} \leq \frac{1}{t^\ep}C(\norm{Z(t)}_{H^k(\Sc)})\norm{Z(t)}_{H^k(\Sc)}.
\end{equation*}. 
\item There exist  constants $m\in \Zbb_{\geq 1}$ and $\delta>0$ such that for any $\mathring{V}=(\mathring{V}^K)\in H^{k+1}(\Sc,\Vbb^N)$ satisfying $\norm{\mathring{V}}_{H^{k+1}(\Sc)}<\delta$, there
exists a unique solution
\begin{equation*} 
V=(V^K) \in C^0\bigl((0,t_0],H^{k+1}(\Sc,\Vbb^N)\bigr)\cap C^1\bigl((0,t_0],H^{k}(\Sc,\Vbb^N)\bigr)
\end{equation*}
to the GIVP \eqref{MwaveJ.1}-\eqref{MwaveJ.2} for the extended system,  where $t_0$ is defined in \eqref{time0}.   The  solution $V$ satisfies the bounds
\begin{gather*}
\norm{V_0(t)}_{L^\infty(\Sc)} \lesssim 1, \quad
\norm{V_0(t)}_{H^k(\Sc)}  \lesssim \frac{1}{t^{\ep}}, \quad 
\norm{\Pbb V(t)}_{H^k(\Sc)}  \lesssim t^{\nu},  \\
\norm{\Dc V(t)}_{H^k(\Sc)}  \lesssim \frac{1}{t^{\kappa}}, \quad
\norm{\Pbb V(t)}_{H^{k-1}(\Sc)}  \lesssim t^{\nu+\kappa-\zc} \AND
\norm{\Dc V(t)}_{H^{k-1}(\Sc)}  \lesssim \frac{1}{t^{\zc}}
\end{gather*}
for $t\in (0,t_0]$.
\item Given  initial data $\mathring{V}$ satisfying  the constraint \eqref{ind1},  then the solution $V$
determines a unique classical solution $\ub^K \in C^2(\Mb_{r_0})$,
with $r_0=\rho_0^m$, of the IVP
\begin{align*}
\gb^{\alpha\beta}\nablab_\alpha \nablab_\beta \ub^K &= \ab^{K\alpha\beta}_{IJ}\nablab_\alpha \ub^I \nablab_\beta \ub^J
\quad \text{in $\Mb_{r_0}$,}\\
(\ub^K, \del{\tb}\ub^K) &= (\vb^K,\wb^K_1) \hspace{1.5cm} \text{in $\Sigmab_{r_0}$,}
\end{align*}
where $\ub^K$, $\vb^K$ and $\wb^K$ are determined from $V$ by \eqref{ubarsl}, \eqref{ind3}-\eqref{indz}, and the solution  $\ub^K$ satisfy the pointwise bounds
\begin{equation*}
|\ub^K|\lesssim \frac{\tb-\rb}{1+\tb-\rb}\biggl(\frac{1}{1+\tb-\rb}\biggr)^{\nu+\kappa-\zc-\frac{1}{2}} \quad \text{in $\Mb_{r_0}$.}
\end{equation*} 
\end{enumerate}
\end{thm}
\begin{proof}
Assuming that the GIVP for the extended system \eqref{MwaveJ.1}-\eqref{MwaveJ.2} satisfies the flow assumptions \eqref{flowassump.1}-\eqref{flowassump2},  and the  constants $\epsilon, \zc, \nu, \kappa$ satisfy the inequalities \eqref{kappa-nu-fix1}, then the  proof follows  directly from Theorem 4.1 from \cite{Oliynyk2021107}.
\end{proof}

\bigskip

\noindent \textit{Acknowledgements:}
The author thanks Todd Oliynyk and Tracey Balehowsky for useful comments and discussions in  the developing of this manuscript.

\appendix
\section{Indexing conventions\label{indexing}}
Below is a summary of the indexing conventions that are employed throughout this article:
\medskip
\begin{center}
\begin{tabular}{|l|c|c|l|} \hline
Alphabet & Examples & Index range & Index quantities  \\ \hline
Lowercase Greek & $\mu,\nu,\gamma$ & $0,1,2,3$ & spacetime coordinate components: $(x^\mu)=(t,r,\theta,\phi)$ \\ \hline
Uppercase Greek &$\Lambda,\Sigma,\Omega$ & $2,3$, & spherical coordinate components:$(x^\Lambda)=(\theta,\phi)$ \\ \hline
Lowercase Latin &$i,j,k $ & $1,2,3$ & spatial coordinate components:$(y^i)=(\rho,\theta,\phi)$ \\ \hline
Uppercase Latin &$I,J,K$ & $1$ to $N$ & wave equation indexing:$u^I$  \\ \hline
Lowercase Calligraphic &$\qc,\pc,\rc$ & 0,1 & time and radial coordinate components:$(x^\qc)=(t,r)$\\ \hline
Uppercase Calligraphic &$\Ic,\Jc,\Kc$ & 0,1,2,3,4 & first order wave formulation indexing:$V^K_\Ic$  \\ \hline
\end{tabular}
\end{center}

\section{Conformal Transformations} \label{ctrans}
In this section, we recall a number of formulas that govern the transformation laws for geometric objects under a conformal transformation that will be needed for our application to wave equations. Under a  
conformal transformation of the form
\begin{equation} \label{gtrans}
\gt_{\mu\nu} = \Omega^2 g_{\mu\nu},
\end{equation}
the Levi-Civita connection $\nablat_\mu$ and $\nabla_\mu$ of $\gt_{\mu\nu}$ and $g_{\mu\nu}$, respectively, are related by
\begin{equation*} 
\nablat_{\mu}\omega_\nu = \nabla_\mu\omega_\nu - \Cc_{\mu\nu}^\lambda \omega_\lambda,
\end{equation*}
where 
\begin{equation*}
\Cc_{\mu\nu}^\lambda = 2\delta^\lambda_{(\mu}\nabla_{\nu)}\ln(\Omega)
-g_{\mu\nu}g^{\lambda\sigma}\nabla_\sigma \ln(\Omega).
\end{equation*}
Using this, it can be shown that the wave operator transforms as
\begin{equation}\label{wavetransA}
\gt^{\mu\nu}\nablat_\mu\nablat_\nu \ut - \frac{n-2}{4(n-1)}\Rt \ut = \Omega^{-1-\frac{n}{2}}
\biggl(g^{\mu\nu}\nabla_\mu\nabla_\nu u - \frac{n-2}{4(n-1)}R u\biggr)
\end{equation}
where $\Rt$ and $R$ are the Ricci curvature scalars of $\gt$ and $g$, respectively, $n$ is the dimension of
spacetime, and
\begin{equation} \label{utransA}
\ut = \Omega^{1-\frac{n}{2}}u.
\end{equation} 
Assuming now that the scalar functions $\ut^K$ satisfy the system of wave equations
\begin{equation} \notag
\gt^{\mu\nu}\nablat_\mu\nablat_\nu \ut^K - \frac{n-2}{4(n-1)}\Rt \ut^K =\ft^K,
\end{equation}
it then follows immediately from \eqref{wavetransA} and
\eqref{utransA} that the scalar functions
\begin{equation} \label{utrans}
u^K = \Omega^{\frac{n}{2}-1}\ut^K
\end{equation}
satisfy the conformal system of wave equations given by
\begin{equation} \notag
g^{\mu\nu}\nabla_\mu\nabla_\nu u^K - \frac{n-2}{4(n-1)} R u^K =f^K
\end{equation}
where
\begin{equation} \label{ftransA}
f^K = \Omega^{1+\frac{n}{2}}\ft^K.
\end{equation}
Specializing to source terms $\ft^K$ that are quadratic in the derivatives, that is, of the form
\begin{equation} \label{Awet3}
\ft^K = \at^{K\mu\nu}_{IJ}\nablat_\mu \ut^I \nablat_\nu \ut^J,
\end{equation}
a short calculation using \eqref{gtrans} and \eqref{utrans} shows that the corresponding conformal source $f^K$,
defined by \eqref{ftransA}, is given by
\begin{align} \label{Nlint}
f^K = \at^{K\mu\nu}_{IJ}\biggl(\Omega^{3-\frac{n}{2}}\nabla_\mu u^I\nabla_\nu u^J& +\bigg(\frac{n}{2}-1\biggr)\Omega^{4-\frac{n}{2}}\bigl(\nabla_\mu\Omega^{-1} u^I \nabla_\nu u^J+ \nabla_\mu u^I \nabla_\nu\Omega^{-1} u^J \bigr) \notag\\
&\qquad  + \bigg(1-\frac{n}{2}\biggr)^2  \Omega^{5-\frac{n}{2}}\nabla_\mu\Omega^{-1} \nabla_\nu\Omega^{-1} u^I u^J\biggr).
\end{align}

\bibliographystyle{unsrt}
\bibliography{refs}

\end{document}